\documentclass[12pt]{article}
\usepackage{epsfig}
\usepackage{amsfonts}
\usepackage{amssymb}
\usepackage{amsmath}
\usepackage{amsthm}
\usepackage{latexsym}
\usepackage{graphicx}
\usepackage{bm}
\usepackage{tabularx}
\usepackage{booktabs}
\usepackage{enumerate}
\usepackage{caption}
\usepackage[titletoc]{appendix}
\usepackage[numbers]{natbib}
\usepackage{bbm}
\usepackage{url}
\usepackage{hyperref}
\usepackage{color, float}

\usepackage{subcaption}
\usepackage{mwe}

\catcode`\@=11

\newcommand{\shortdot}[1]{\raisebox{-0.4pt}{$\stackrel{\bullet}{#1}$}}
\newcommand{\updot}[1]{\raisebox{0.9pt}{$\stackrel{\bullet}{#1}$}} 

\newcommand{\paren}[1]{\left(#1\right)}
\newcommand{\sqparen}[1]{\left[#1\right]}

\newcommand{\E}[1]{{\mathrm{E}\left[#1\right]}}

\newtheorem{theorem}{Theorem}[section]
\newtheorem{lemma}[theorem]{Lemma}
\newtheorem{corollary}[theorem]{Corollary}
\newtheorem{proposition}[theorem]{Proposition}

\theoremstyle{remark}
\newtheorem{rem}[theorem]{Remark}

\title{New Perspectives on the Erlang-A Queue}

\author{
Andrew Daw \\ School of Operations Research and Information Engineering \\ Cornell University
\\ 257 Rhodes Hall, Ithaca, NY 14853 \\  amd399@cornell.edu \\
 \and
  Jamol Pender \\ School of Operations Research and Information Engineering \\ Cornell University
\\ 228 Rhodes Hall, Ithaca, NY 14853 \\  jjp274@cornell.edu
 }

\begin{document}

\maketitle

\begin{abstract}

The non-stationary Erlang-A queue is a fundamental queueing model that is used to describe the dynamic behavior of large scale multi-server service systems that may experience customer abandonments, such as call centers, hospitals, and urban mobility systems.  In this paper, we develop novel approximations to all of its transient and steady state moments, the moment generating function, and the cumulant generating function.  We also provide precise bounds for the difference of our approximations and the true model.  More importantly, we show that our approximations have \emph{explicit stochastic representations as shifted Poisson random variables}.  Moreover, we are also able to show that our approximations and bounds also hold for non-stationary Erlang-B and Erlang-C queueing models under certain stability conditions.  Finally, we perform numerous simulations to support the conclusions of our results.
\end{abstract}




\section{Introduction}

Markov processes are important modeling tools that help researchers describe real-world phenomena.  Thus, it comes as no surprise that the Erlang-A model, which is a Markovian and multi-server queueing model that incorporates customer abandonments, is an important modeling tool in a multitude of application settings.  Some of the more prominent applications include telecommunications, healthcare, urban mobility and transportation, and more recently cloud computing.  See for example the following work by \citet{mandelbaum2002queue, massey2002analysis, yom2014erlang, pender2016law}.  Despite its importance in many different applications, the Erlang-A queueing model has remained to be  very difficult to analyze and understand.  Even the analysis of the moments of Erlang-A queue beyond the fourth moment has remained an important topic for additional study.

It is well known that the stationary setting of the Erlang-A is much easier to analyze than its non-stationary counterpart. Some common approaches used to analyze non-stationary and state dependent queueing models including asymptotic methods such as heavy traffic limit theory and strong approximations theory, see for example \citet{halfin1981heavy, mandelbaum1998strong}.  Uniform acceleration is extremely useful for approximating the transition probabilities and moments such as the mean and variance of Markov processes.  Moreover, the strong approximation methods are useful for analyzing the sample path behavior of the Markov process by showing that the sample paths of properly rescaled queueing processes converge to deterministic dynamical systems and Gaussian process limits.

However, there are two main drawbacks of these asymptotic methods.  The first is that the method is asymptotic as a function of the model parameters and the results really only hold when the rates are large and are nearly infinite.  Thus, the quality of the approximations depends significantly on the size of the model parameters and these asymptotic methods have been shown to be quite inaccruate for moderate sized model parameter settings, see for example \citet{massey2011poster, massey2013gaussian}.  The second main drawback is that the asymptotic methods do not generate any important insights for the moments or cumulant moments beyond order two since the limits are are based on Brownian motion.  Since Brownian motion has symmetry, its cumulants are all zero beyond the second order.  Thus, Brownian approximations are limited in their power to capture asymmetries in higher moments or even the dynamics of the moment generating function, cumulant generating function, or Fourier transform.  Moreover, it has been shown recently by \citet{pender2014gram, engblom2014approximations} that the Erlang-A and its variants have non-trivial amounts of skewness and excess kurtosis, which implies that the Erlang-A are not nearly Gaussian for moderate sized queues.  These results also demonstrate that it is important to capture the behavior of the Erlang-A model beyond its second moment as this information can be used in staffing decisions \citet{massey2017performance}.

One common approximation method that is used in the stochastic networks, queueing, and chemical reactions literature is a \emph{moment closure approximation}.  Moment closure approximations are used to approximate the moments of the queueing process with a surrogate distribution.  It is often the case that the set of moment equations for a large number of queueing models are not closed, see for example \citet{matis2001transient, pender2014laguerre}.   Thus, the closure approximation helps approximate the moments with a closed system using the surrogate distribution.  One such method used by \citet{ pender2016sampling,  pender2017approximations} is to use Hermite polynomials for approximating the distribution of the queue length process. In fact, they show that using a quadratic polynomial works quite well.  Since the Hermite polynomials are orthogonal to the Gaussian distribution, which has support on the entire real line, these Hermite polynomial approximations do not take into account the discreteness of the queueing process and the fact the queueing process is non-negative.  However, they show that Hermite polynomials are natural to analyze since they are orthogonal with respect to the Gaussian distribution and the heavy traffic limits of multi-server queues are Gaussian.

In this paper, we perform an in-depth analysis of the moments and the moment generating function of the non-stationary Erlang-A queue.  As the Erlang-B and Erlang-C queueing models are special cases of the Erlang-A model, we are able to obtain similar results for those models.  Our approach is to use convexity and exploit Jensen's and the FKG inequality to obtain bounds on the moments and moment generating function of the Erlang-A queue.  What we find even more exciting is that we are able to provide a stochastic representation of our approximations and bounds as a Poisson random variables with a constant shift.  This shifted Poisson was observed in peer to peer networks by \citet{ferragut2012content}, however, we will show in the sequel, this novel representation will allow us to view our bounds and approximations in a new way.


 \subsection{Main Contributions of the Paper}

The main contributions of this work can be summarized as follows:
\begin{itemize}
\item We provide new approximations for the moments, moment generating function, and cumulant generating function for the nonstationary Erlang-A queue exploting FKG and Jensen's inequalities.
\item We derive a novel stochastic intepretation and representation of our approximations as shifted Poisson random variables or $M/M/\infty$ queues, depending on the context.  This sheds new light on the complexity of queues in heavy traffic or critically loaded regimes.
\item We prove precise error bounds for our approximations and we also prove new upper and lower bounds for the nonstationary Erlang-A queue that become exact in certain parameter settings.
\end{itemize}


\subsection{Organization of the Paper}

The remainder of this paper is organized as follows. Section~\ref{secQMod} introduces the nonstationary Erlang-A queueing model and its importance in stochastic network theory.  In Section~\ref{MeanJensen}, we provide approximations for the moments of the Erlang-A system and use these to bound the true values. In Section~\ref{MGFsec} we derive approximations for the moment generating function and cumulant moment generating function of the Erlang-A queue. We again bound the true values by these approximations, and we also find a representation for our approximations in terms of Poisson random variables or $M/M/\infty$ queues, depending on the context.

\section{The Erlang-A Queueing Model}\label{modelsec}

\label{secQMod}

The Erlang-A queueing model is a fundamental queueing model in the stochastic processes literature.  The work of \citet{mandelbaum1998strong}, shows that the $M(t)/M/c+M$ queueing system process $ Q\equiv \{ Q(t) | t \geq 0 \} $ is represented by the following stochastic, time changed integral equation:
$$
Q(t) = Q(0) + \mathit{\Pi}_1 \left(\int^{t}_{0} \lambda(s) ds \right) -  \mathit{\Pi}_2 \left(\int^{t}_{0} \mu \cdot (Q(s) \wedge c )ds \right)
- \mathit{\Pi}_3 \left(\int^{t}_{0} \theta \cdot (Q(s) - c )^+ds \right) ,
$$
where $  \mathit{\Pi}_i \equiv \{  \mathit{\Pi}_i(t) | t \geq 0 \} $
for $ i = {1,2,3} $ are i.i.d.\  standard (rate 1) Poisson processes.  Thus, we can write the sample path dynamics of the Erlang-A queueing process in terms of three independent unit rate Poisson processes.  A deterministic time change for $\mathit{\Pi}_1 $ transforms it into a non-homogeneous Poisson arrival process with rate $ \lambda(t)$ that counts the customer arrivals that occured in the time interval [0,t).  A random time change for the Poisson process $\mathit{\Pi}_2$ , gives us a departure process that counts the number of serviced customers.  We implicitly assume that the number of servers is $c \in \mathbb{Z}^+$ and that each server works at rate $\mu$. Finally, a the random time change of $ \mathit{\Pi}_3$  gives us a counting process for the number of customers that abandon service.  We also assume that the abandonment distribution is exponential and the rate of abandonments is equal to $\theta$.


One of the main reasons that the Erlang-A queueing model has been studied so extensively is because several important queueing models are special cases of it.  One special case is the infinite server queue.  The infinite server queue can be derived from the Erlang-A queue in two ways.  The first way is to set the number of servers to infinity.  This precludes any abandonments since the abandonment rate $\theta \cdot (Q(t) - c)^+$ is always equal to zero when the number of servers is infinite.  The second way to derive the infinite server queue is to set the service rate $\mu$ equal to the abandonment rate $\theta$.  When $\mu = \theta$, this implies that the sum of the service and abandonment departure processes is equal to a linear function i.e. $\mu \cdot (Q(t) \wedge c) + \theta \cdot (Q(t) - c )^+ = \mu \cdot Q(t) = \theta \cdot Q(t)$.  Thus, the Erlang-A queueing model becomes an infinite server queue.

One of the main and important insights of \citet{halfin1981heavy} is that for multi-server queueing systems, it is natural to scale up the arrival rate and the number of servers simultaneously.  This scaling known as the \emph{Halfin-Whitt} scaling and been an important modeling technique for modeling call centers in the queueing literature.  Since the $M(t)/M/c+M$ queueing process is a special case of a single node \emph{Markovian service network},
we can also construct an associated, \emph{uniformly accelerated} queueing process where both the new arrival rate
$\eta\cdot\lambda(t)$ and the new number of servers $\eta\cdot c$ are both scaled by the same factor $\eta>0$.  Thus, using the \emph{Halfin-Whitt} scaling for the Erlang-A model, we arrive at the following sample path representation for the queue length process as
 \begin{eqnarray*}
Q^{\eta}(t) &=& Q^{\eta}(0) + \mathit{\Pi}_1 \left(\int^{t}_{0}  \eta \cdot \lambda(s) ds \right) -  \mathit{\Pi}_2 \left(\int^{t}_{0} \mu \cdot (Q^{\eta}(s) \wedge  \eta \cdot c )ds \right) \\
&&- \mathit{\Pi}_3 \left(\int^{t}_{0} \theta \cdot (Q^{\eta}(s) - \eta \cdot c )^+ds \right) \\
&=& Q^{\eta}(0) + \mathit{\Pi}_1 \left( \int^{t}_{0}  \eta \cdot \lambda(s) ds \right) -  \mathit{\Pi}_2 \left(\int^{t}_{0}  \eta \cdot \mu \cdot \left( \frac{Q^{\eta}(s)}{\eta} \wedge c \right)ds \right) \\
&&- \mathit{\Pi}_3 \left(\int^{t}_{0} \eta \cdot \theta \cdot \left( \frac{Q^{\eta}(s)}{\eta} - c \right)^+ds \right) .
\end{eqnarray*}

The \emph{Halfin-Whitt} scaling is defined by  simultaneously scaling up the rate of customer demand (which is the arrival rate) with the number of servers.  In the context of call centers this is scaling up the number of customers and scaling up the number of agents to answer the phones.  In the context of hospitals or healthcare this might be scaling up the number of patients with the number of beds or nurses.  Taking the following limits gives us the \emph{fluid} models of \citet{mandelbaum1998strong}, i.e.
 \begin{equation}
\lim_{\eta\to\infty} \frac{1}{\eta} Q^{\eta}(t) =   q(t) \hspace{3mm}
\mathrm{a.s.}
\end{equation}
where the deterministic process $q(t)$, the \emph{fluid mean}, is
governed by the one dimensional ordinary differential equation (ODE)
\begin{equation}
\label{fldmean}
 \shortdot{q}(t) = \lambda(t) - \mu \cdot (q(t) \wedge c) - \theta \cdot (q(t) - c)^+ .
\end{equation}
Moreover, if one takes a diffusion limit i.e.
 \begin{equation}
\lim_{\eta\to\infty} \sqrt{\eta} \left( \frac{1}{\eta} Q^{\eta}(t) -   q(t) \right) \Rightarrow \tilde{Q}(t) \hspace{3mm}
\end{equation}
one gets a diffusion process where the variance of the diffusion is given by the following ODE
\begin{eqnarray}
\label{diffvar}
 \updot{\mathrm{Var}}\left[ \tilde{Q}(t) \right]
 & =& \lambda(t) + \mu \cdot (q(t) \wedge c) + \theta \cdot (q(t) - c)^+  \nonumber \\
&&- 2 \cdot  \mathrm{Var}\left[\tilde{Q}(t) \right] \cdot \left( \mu \cdot \{ q(t) < c \} +   \theta \cdot \{ q(t) \geq  c \} \right) .
\end{eqnarray}

\subsection{Mean Field Approximation is Identical to the Fluid Limit}

\begin{figure}[H]
\centering
\vspace{-1in}
\hspace{-.25in}~\includegraphics[scale = .4]{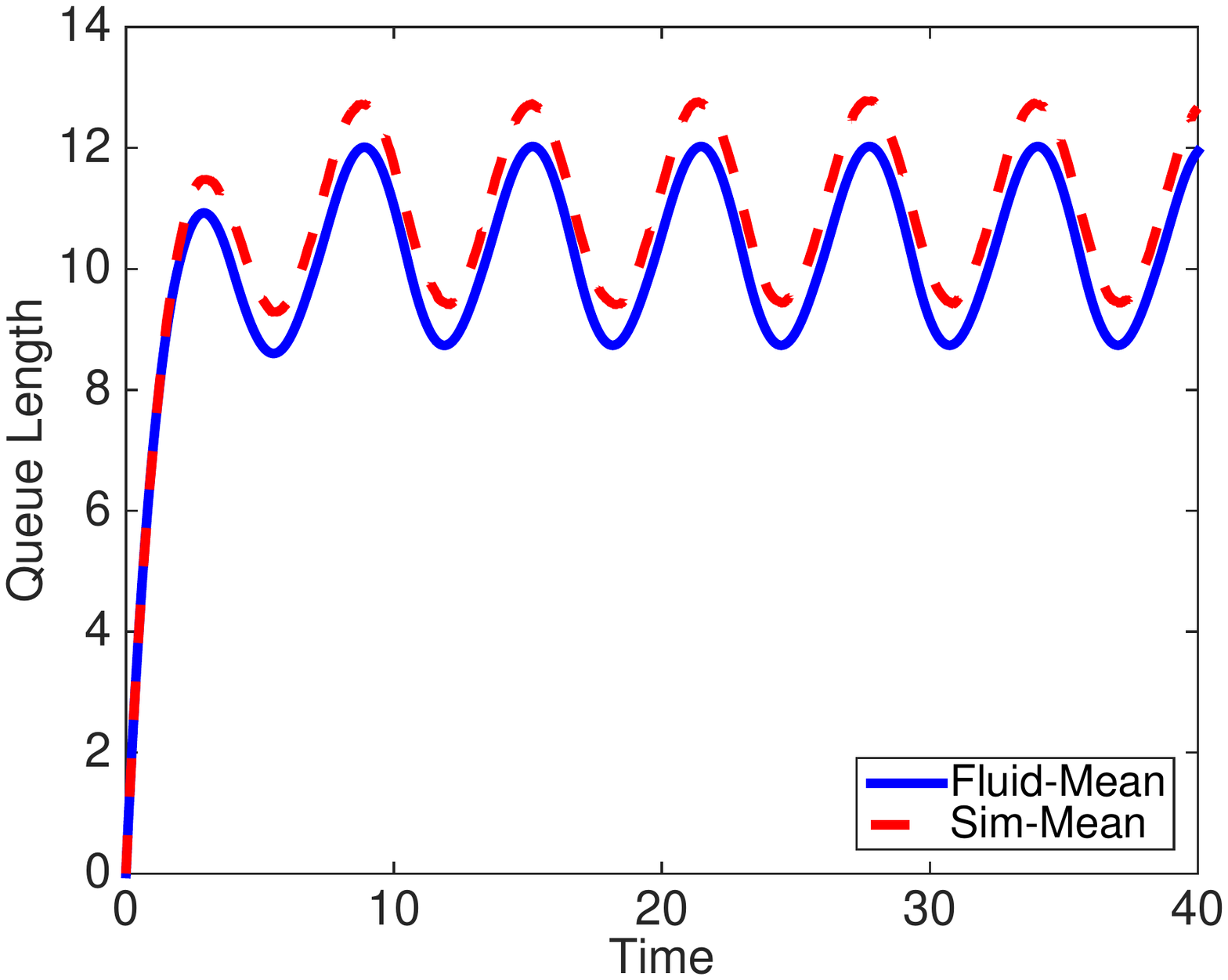}~\hspace{-.15in}~\includegraphics[scale = .4]{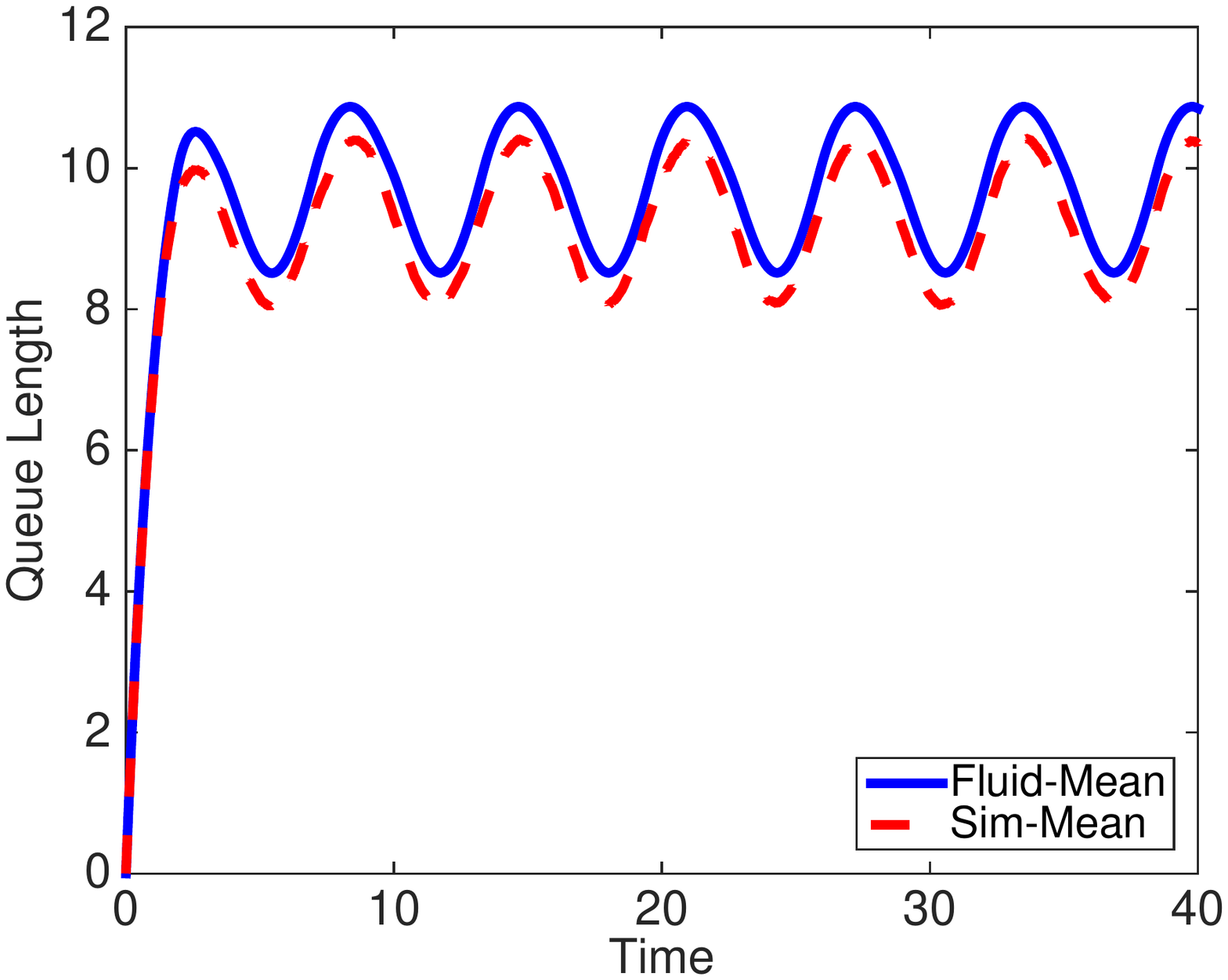}
\captionsetup{justification=centering}
\vspace{-1.25in}
 \caption{$\lambda(t) = 10 + 2 \cdot \sin (t) $, $\mu = 1$,  $Q(0) = 0$, $c=10$. \\
 $\theta = 0.5$ (Left)  and $\theta = 2 $ (Right). \label{Fig1} }
\end{figure}

In addtion to using strong approximations to analyze the queue length process one can also use the functional Kolmogorov forward equations as outlined in \citet{massey2013gaussian}.  The functional forward equations for the Erlang-A model are derived as,
\begin{eqnarray}
\label{FOREQN}
\updot{\mathrm{E}} [f(Q(t)) ] &\equiv& \frac{d}{dt} \mathrm{E}[ f(Q(t)) | Q(0) = q(0) ] \\
&=& \lambda \cdot
\mathrm{E}\sqparen{f(Q(t)+1)-f(Q(t))}
+
\mathrm{E}\sqparen{\delta\paren{Q(t), c}\cdot\paren{f(Q(t)-1)-f(Q(t))}} ,
\end{eqnarray}
for all appropriate functions $f$ and where $ \delta\paren{Q(t), c} = \mu \cdot (Q(t) \wedge c) + \theta \cdot (Q(t) - c)^+ $.  For the special case where $f(x) = x$, we can derive an ode for the mean queue length process as
  \begin{eqnarray}\label{meanforeqn}
     \updot{\mathrm{E}} [Q(t)]  &=&\lambda(t) - \mu \cdot \mathrm{E}[ ( Q(t) \wedge c ) ]  -  \theta \cdot \mathrm{E}[ ( Q(t) - c )^+ ]  .
  \end{eqnarray}
  The first thing to note is that this equation is not autonomous and one needs to know the distribution of $Q(t)$ a priori in order to compute the expectations on the righthand side of Equation \ref{meanforeqn}.  To know the distribution a priori is impossible except in some special cases like the infinite server setting.  However, it is easy to derive simple approximations for the mean queue length by making some assumptions on the queue length process.  This is known as a closure approximation and one common closure approximation method is to simply take the expectations from outside the function to inside the function.  This implies that the expectation $E[f(X)]$ becomes $f(E[X])$.  This method is known as a mean field approximation in physics and is also known as the deterministic mean approximation of \citet{massey2013gaussian}.   By applying the mean field approximation to Equation \ref{meanforeqn}, we can show that the resulting differential equation is given by the following autonomous ODE
  \begin{eqnarray} \label{fluidmean}
     \updot{\mathrm{E}} [Q_{f}(t)]  &=&\lambda(t) - \mu \cdot ( \mathrm{E}[Q_{f}] \wedge c )   -  \theta \cdot  (\mathrm{E}[Q_{f}] - c )^+  .
  \end{eqnarray}
By careful inspection, one can observe that the ode given by the mean field approximation is identical to the fluid limit of \ref{fldmean}.   Moreover, if one simulates the queueing process and compare it to the mean field limit, one notices an ordering property.  For example on the left of Figure \ref{Fig1}, we simulate the Erlang-A queue and compare to the fluid model.  We observe that when $\theta < \mu$, that the simulated mean is larger than the fluid mean.  This is precisely what our results predict.  Moreover, on the right of Figure \ref{Fig1}, we simulate the Erlang-A queue and compare to the fluid model when $\theta > \mu$ and observe that simulated queue length is smaller than the fluid limit.

Our goal in this work is to explain the behavior that we observe in Figure~\ref{Fig1}, which we will do in the following section. Before concluding our overview of the Erlang-A queueing model, we make a brief remark for notational clarity.

\begin{rem}\label{notationremark}
Throughout the remainder of this work, we use $Q(t)$ to represent the true queueing process and $Q_f(t)$ to represent the fluid approximation of it. This fluid approximation is a stochastic process that will be fully described in this work. In fact, in Section~\ref{MGFsec} we use characterize the fluid approximations and use insight from these representations to bound the true queue length from above and below.
\end{rem}

\section{Inequalities for the Moments of the Erlang-A Queue} \label{MeanJensen}

In this section, we prove when the true moments of the Erlang-A queue are either dominated or dominates their corresponding fluid limit.  We find that the relationship between the service rate and the abandonment rate determines whether or not the moment is dominated by the fluid limit. This section is organized as follows. In Subsection~\ref{meanineq}, we derive inequalities for the true mean of the Erlang-A and its fluid approximation. In Subsection~\ref{MomentJensen} we extend these inequalities to analogous results for the $m^\text{th}$ moment of the queueing system. Finally, in Subsection~\ref{momentnumerical} we provide figures from numerical experiments that demonstrate these findings.

\subsection{Inequalities for the Mean}\label{meanineq}

  We begin with analysis of the mean of the Erlang-A queue. Before we proceed, we first establish a lemma for comparisons of ordinary differential equations that will be fundamental to our approach to the results.

\begin{lemma}[A Comparison Lemma]
Let $f: \mathbb{R}^2 \to \mathbb{R}$ be a continuous function in both variables.  If we assume that initial value problem
\begin{equation}
\shortdot{x}(t) = f(t, x(t)), \ x(0) = x_0
\end{equation}
has a unique solution for the time interval [0,T] and
\begin{equation}
\shortdot{y}(t) \leq f(t, y(t)) \quad \mathrm{for} \ t \in [0,T] \ \mathrm{and} \ y(0) \leq x_0
\end{equation}
then $ x(t) \geq y(t) $ for all $t \in [0,T]$.
\begin{proof}
The the proof of this result is given in \citet{hale2013introduction}.
\end{proof}
\end{lemma}

With this lemma in hand, we can now derive relationships for the fluid limit and the true mean. As seen in the proof, these results follow from the application of this differential equation comparison lemma and the convexity seen in the fluid approximation.

\begin{theorem}\label{jensen1}
For the Erlang-A queue, if $Q(0) = Q_{f}(0)$, then the true mean dominates the fluid limit when $\theta < \mu$, the fluid limit dominates the true mean when $\theta > \mu$, and the two means are equal when $\theta = \mu$.

\begin{proof}
Recall that the true mean satisfies the following differential equation
  \begin{eqnarray*}
     \updot{\mathrm{E}} [Q(t)]  &=&\lambda(t) - \mu \cdot \mathrm{E}[ ( Q \wedge c ) ]  -  \theta \cdot \mathrm{E}[ ( Q - c )^+ ]
  \end{eqnarray*}
  and the fluid limit satisfies the following differential equation
    \begin{eqnarray*}
     \updot{\mathrm{E}} [Q_{f}(t)]  &=& \lambda(t) - \mu \cdot ( \mathrm{E}[Q_{f}] \wedge c )   -  \theta \cdot  (\mathrm{E}[Q_{f}] - c )^+  .
  \end{eqnarray*}
  We can simplify both equations by observing that $( X \wedge c ) + ( X - c )^+ = X $ for any random variable X.  Thus, we have the following two equations for the true mean and the fluid limit
    \begin{eqnarray*}
     \updot{\mathrm{E}} [Q(t)]  &=&\lambda(t) - \theta \cdot E[ Q]  +  (\theta - \mu) \cdot \mathrm{E}[ ( Q \wedge c ) ]
\\
     \updot{\mathrm{E}} [Q_{f}(t)]  &=& \lambda(t) - \theta \cdot \mathrm{E}[Q_{f}]   + (\theta - \mu) \cdot  (\mathrm{E}[Q_{f}] \wedge c ) .
  \end{eqnarray*}

  If we take the difference of the two equations, we obtain the following
    \begin{eqnarray*}
     \updot{\mathrm{E}} [Q(t)] - \updot{\mathrm{E}} [Q_{f}(t)]  &=& \lambda(t) - \theta \cdot E[ Q]  +  (\theta - \mu) \cdot \mathrm{E}[ ( Q \wedge c ) ]    \\
     &-& \lambda(t) + \theta \cdot \mathrm{E}[Q_{f}]    - (\theta - \mu) \cdot  (\mathrm{E}[Q_{f}] \wedge c )  \\
     &=&   \theta \cdot \left(   \mathrm{E}[Q_{f}]   - E[ Q ] \right)  + ( \theta - \mu) \cdot  \left( \mathrm{E}[ ( Q \wedge c ) ]    - ( \mathrm{E}[Q_{f}]\wedge c ) \right) \\
  \end{eqnarray*}
  Now since the minimum function $(Q \wedge c)$ is a concave function, we have that
      \begin{eqnarray*}
       \left( \mathrm{E}[ ( Q \wedge c ) ]    - ( \mathrm{E}[Q] \wedge c ) \right) &\leq& 0
  \end{eqnarray*}
  for any random variable Q.  Thus, we have that for $\theta < \mu$
      \begin{eqnarray*}
     \updot{\mathrm{E}} [Q(t)] - \updot{\mathrm{E}} [Q_{f}(t)]  &\geq& 0,
  \end{eqnarray*}
 and for $\theta > \mu$
        \begin{eqnarray*}
     \updot{\mathrm{E}} [Q(t)] - \updot{\mathrm{E}} [Q_{f}(t)]  &\leq& 0.
  \end{eqnarray*}
Finally, for $\theta = \mu$, we have that
        \begin{eqnarray*}
     \updot{\mathrm{E}} [Q(t)] - \updot{\mathrm{E}} [Q_{f}(t)]  &=& 0
  \end{eqnarray*}
  since both differential equations are initialized with the same value and the origin is an equilibrium point for the difference.  This completes the proof.
\end{proof}
\end{theorem}

As discussed in Section~\ref{modelsec}, the Erlang-A model is quite versatile in its relation to other queueing systems of practical interest. In the two following corollaries, we find that Theorem~\ref{jensen1} can be applied to the Erlang-B and Erlang-C models.

\begin{corollary}
For the Erlang-B queueing model,  if $Q(0) = Q_{f}(0)$, then $\mathrm{E}[Q(t)] \leq  \mathrm{E}\left[Q_{f}(t) \right]$ for all $t \geq 0$.

\begin{proof}
This is obvious after noticing that the Erlang-B queue is a limit of the Erlang-A queue by letting $\theta \to \infty$.
\end{proof}
\end{corollary}

\begin{corollary}
For the Erlang-C queueing model,  if $Q(0) = Q_{f}(0)$, then $\mathrm{E}[Q(t)] \geq  \mathrm{E}\left[Q_{f}(t) \right]$ for all $t \geq 0$.

\begin{proof}
This is obvious after noticing that the Erlang-C queue is an Erlang-A queue with $\theta=0$.  Since $\mu $ is assumed to be positive, then we fall into the case where $\theta < \mu$ and this completes the proof.
\end{proof}
\end{corollary}

\begin{rem}
Given that we use Jensen's inequality and the FKG inequality later on in the paper, we find it important to differentiate them.  Here we give an example that sets the two apart.  If we have the following function $Q^n$, then Jensen's inequality implies that $\mathbb{E}[Q^n] \geq \mathbb{E}[Q]^n $.  However, FKG implies that $\mathbb{E}[Q^n] \geq  \mathbb{E}[Q^{n-1}] \cdot \mathbb{E}[Q] $.  We find it interesting that by iterating the FKG inequality $n-2$ more times, it yields Jensen's inequality for the moments of random variables.
\end{rem}



\subsection{Inequalities for the $m^{th}$ Moment}\label{MomentJensen}

In this subsection we will now extend the previous findings for the mean to higher moments of the queueing system. Like the result for the mean, this is again built through observation of the convexity in the differential equation of the fluid approximation.

\begin{theorem}\label{mMoment}
For the Erlang-A queue and $m \in \mathbb{Z}^+$, if $Q(0) = Q_f(0)$, then $\E{Q^m(t)} \geq \E{Q_f^m(t)}$ when $\theta < \mu$, $\E{Q^m(t)} \leq \E{Q_f^m(t)}$ when $\theta > \mu$, and $\E{Q^m(t)} = \E{Q_f^m(t)}$ when $\theta = \mu$.
\begin{proof}
We will use proof by induction. For the base case we can apply Theorem \ref{jensen1}. Now, suppose that the statement holds for $j \in \{1,2, \dots, m -1\}$.  Recall that the $m^\text{th}$ moment satisfies
\begin{align*}
\updot{ \mathrm{E}} \left[Q^m(t)\right]
&=
\lambda(t) \E{\sum_{j=0}^m {m \choose j} Q^j(t) - Q^m(t)}
\\&\quad
+
\E{\left(\sum_{j=0}^m {m \choose j} (-1)^{m-j}Q^j(t) - Q^m(t)\right)\big(\theta Q(t) - (\theta - \mu)(Q(t) \wedge c)\big)}
\\
&=
\lambda(t) \sum_{j=0}^{m-1} {m \choose j} \E{Q^j(t)}
+
\theta  \sum_{j=0}^{m-1} {m \choose j} (-1)^{m-j} \E{Q^{j+1}(t)}
\\&\quad
 +
 (\theta - \mu)\E{\left(\sum_{j=0}^{m-1} {m \choose j} (-1)^{m-1-j}Q^{j+1}(t) \right) \wedge \left(c\sum_{j=0}^{m-1} {m \choose j} (-1)^{m-1-j}Q^j(t) \right)}
\end{align*}
and the approximate autonomous version satisfies
\begin{align*}
\updot{ \mathrm{E}} \left[Q_f^m(t)\right]
&=
\lambda(t) \sum_{j=0}^{m-1} {m \choose j} \E{Q_f^j(t)}
+
\theta  \sum_{j=0}^{m-1} {m \choose j} (-1)^{m-j} \E{Q_f^{j+1}(t)}
\\&\quad
 +
 (\theta - \mu)\sum_{j=0}^{m-1} {m \choose j} (-1)^{m-1-j}\left(\E{Q^{j+1}(t) } \wedge \E{c Q^j(t)} \right)
 \\
&=
\lambda(t) \sum_{j=0}^{m-1} {m \choose j} \E{Q_f^j(t)}
+
\theta  \sum_{j=0}^{m-1} {m \choose j} (-1)^{m-j} \E{Q_f^{j+1}(t)}
\\&\quad
 +
 (\theta - \mu)\left(\E{\sum_{j=0}^{m-1} {m \choose j} (-1)^{m-1-j}Q^{j+1}(t) } \wedge \E{c\sum_{j=0}^{m-1} {m \choose j} (-1)^{m-1-j}Q^j(t)} \right)
\end{align*}
Now by taking the difference, we have that
\begin{align*}
\updot{ \mathrm{E}} \left[Q^m(t)\right] - \updot{ \mathrm{E}} \left[Q_f^m(t)\right]
&=
\lambda(t) \sum_{j=0}^{m-1} {m \choose j} \E{Q^j(t) - Q_f^j(t)}
+
\theta  \sum_{j=0}^{m-1} {m \choose j} (-1)^{m-j} \E{Q^{j+1}(t) - Q_f^{j+1}(t)}
\\&\,\,
 +
 (\theta - \mu)
 \Bigg(
 \E{\left(\sum_{j=0}^{m-1} {m \choose j} (-1)^{m-1-j}Q^{j+1}(t) \right) \wedge \left(c\sum_{j=0}^{m-1} {m \choose j} (-1)^{m-1-j}Q^j(t) \right)}
 \\&
 -
 \E{\sum_{j=0}^{m-1} {m \choose j} (-1)^{m-1-j}Q^{j+1}(t) } \wedge \E{c\sum_{j=0}^{m-1} {m \choose j} (-1)^{m-1-j}Q^j(t)}
 \Bigg).
\end{align*}
Because the minimum is a concave function, we have that for any $X$ and $Y$ with real means $\E{X \wedge Y} \leq \E{X} \wedge \E{Y}$. Thus, we have that for $\theta > \mu$,
$$
\updot{ \mathrm{E}} \left[Q^m(t)\right] - \updot{ \mathrm{E}} \left[Q_f^m(t)\right] \geq 0,
$$
if $\theta < \mu$,
$$
\updot{ \mathrm{E}} \left[Q^m(t)\right] - \updot{ \mathrm{E}} \left[Q_f^m(t)\right] \leq 0,
$$
and if $\theta = \mu$,
$$
\updot{ \mathrm{E}} \left[Q^m(t)\right] = \updot{ \mathrm{E}} \left[Q_f^m(t)\right] = 0
$$
since both differential equations are initialized with the same value, the origin is an equilibrium point for the difference, and all the lower-power terms in the differential equations follow this structure, which we know from the inductive hypothesis. Therefore we see this holds for $m$, which completes the proof.
\end{proof}
\end{theorem}

Again as we have seen for the mean, we can exploit the versatility of the Erlang-A queue to extend these insights to the Erlang-B and Erlang-C models as well.

\begin{corollary}
For the Erlang-B queueing model, if $Q(0) = Q_{f}(0)$, then $\mathrm{E}[Q^m(t)] \leq  \mathrm{E}\left[Q_{f}^m(t) \right]$ for all $t \geq 0$ and $m \in \mathbb{Z}^+$.

\begin{proof}
This is obvious after noticing that the Erlang-B queue is a limit of the Erlang-A queue by letting $\theta \to \infty$.
\end{proof}
\end{corollary}

\begin{corollary}
For the Erlang-C queueing model, if $Q(0) = Q_{f}(0)$, then $\mathrm{E}[Q^m(t)] \geq  \mathrm{E}\left[Q_{f}^m(t) \right]$ for all $t \geq 0$ and $m \in \mathbb{Z}^+$.

\begin{proof}
This is obvious after noticing that the Erlang-C queue is an Erlang-A queue with $\theta=0$.  Since $\mu $ is assumed to be positive, then we fall into the case where $\theta < \mu$ and this completes the proof.
\end{proof}
\end{corollary}

\subsection{Numerical Results} \label{momentnumerical}

\begin{figure}[H]
        \centering
          \begin{subfigure}[b]{0.475\textwidth}
            \centering   \vspace{-.75in}
           \includegraphics[scale = .3]{./Figure_1.pdf}
           \vspace{-.8in}
            \caption{{\small First Moment }}
            \label{fig:1}
        \end{subfigure}
        \hfill
        \begin{subfigure}[b]{0.475\textwidth}
            \centering \vspace{-.75in}
            \includegraphics[scale = .3]{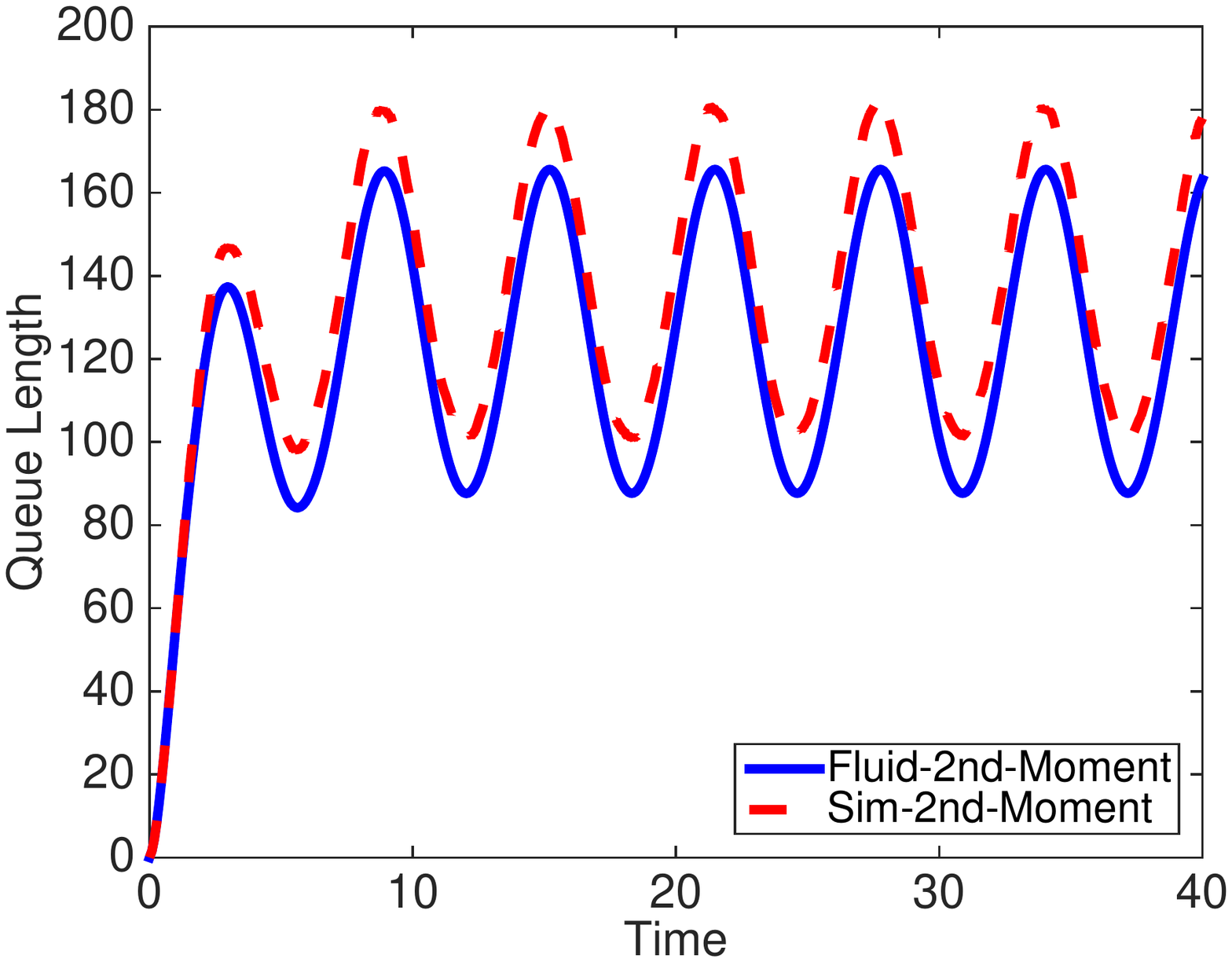}
              \vspace{-.8in}
            \caption{{\small Second Moment }}
            \label{fig:2}
        \end{subfigure}
        \begin{subfigure}[b]{0.475\textwidth}
            \centering   \vspace{-.75in}
            \includegraphics[scale = .3]{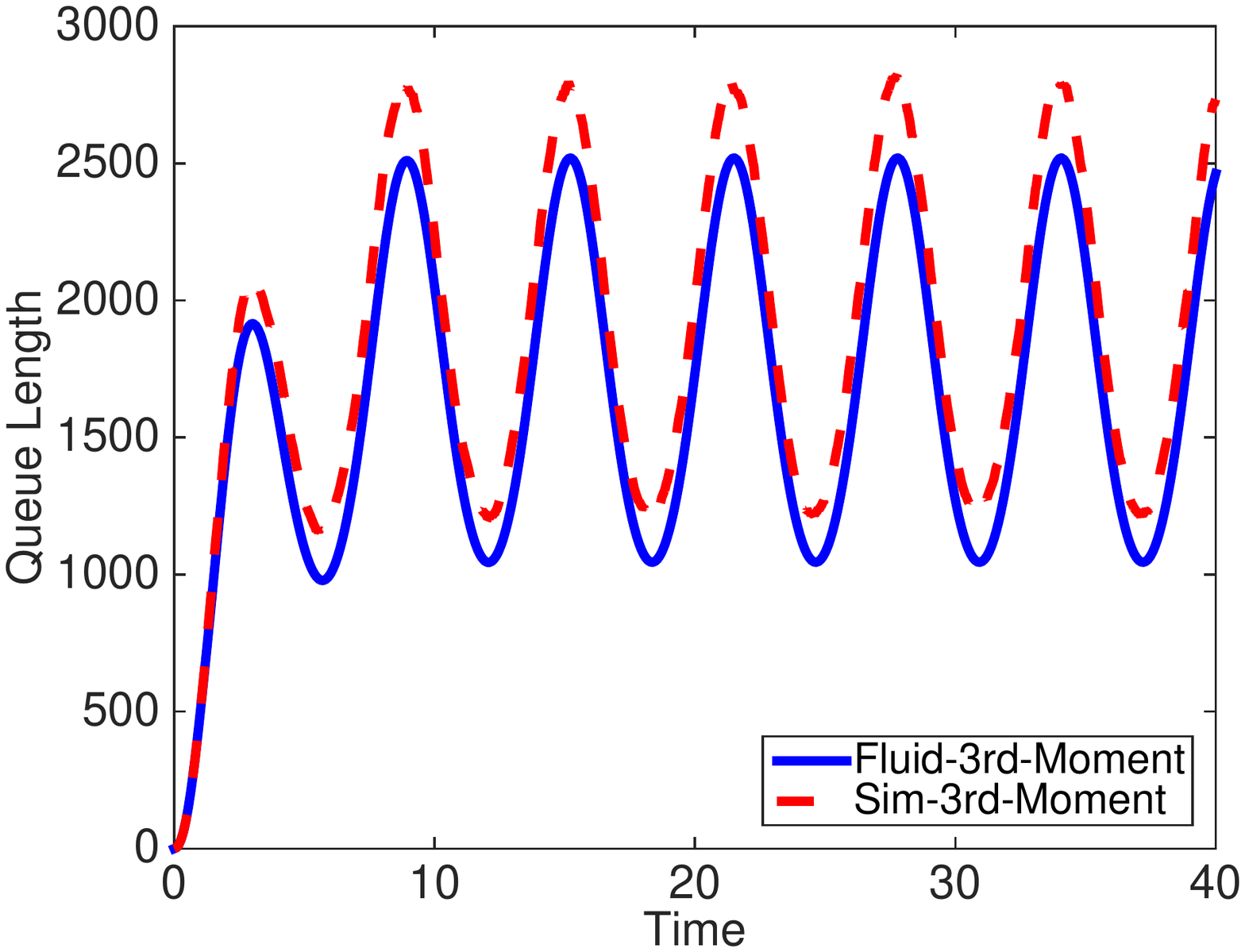}
              \vspace{-.8in}
            \caption{{\small Third Moment }}
            \label{fig:3}
        \end{subfigure}
        \hfill
        \begin{subfigure}[b]{0.475\textwidth}
            \centering   \vspace{-.75in}
            \includegraphics[scale = .3]{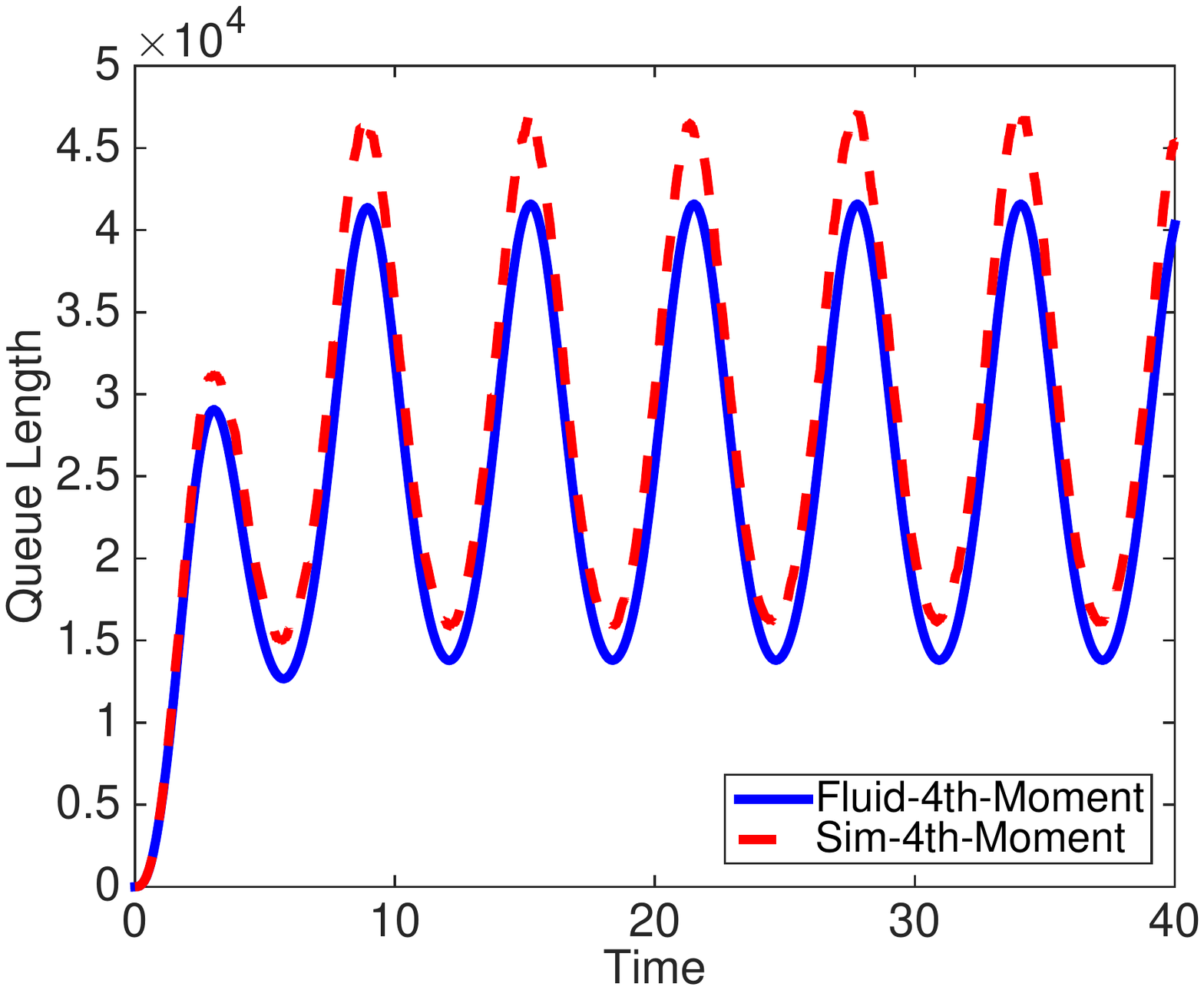}
              \vspace{-.8in}
            \caption{{\small Fourth Moment }}
            \label{fig:4}
        \end{subfigure}
        \caption
        { $\lambda(t) = 10 + 2 \cdot \sin( t) $, $\mu = 1$, $\theta = 0.5$, $Q(0) = 0$, $c=10$.
 }
        \label{fig:1to9pt1}
     \end{figure}

\begin{figure}[H]
        \begin{subfigure}[b]{0.475\textwidth}
            \centering   \vspace{-.75in}
           \includegraphics[scale = .3]{./Figure_6.pdf}
           \vspace{-.8in}
            \caption{{\small First Moment }}
            \label{fig:1}
        \end{subfigure}
        \hfill
        \begin{subfigure}[b]{0.475\textwidth}
            \centering \vspace{-.75in}
            \includegraphics[scale = .3]{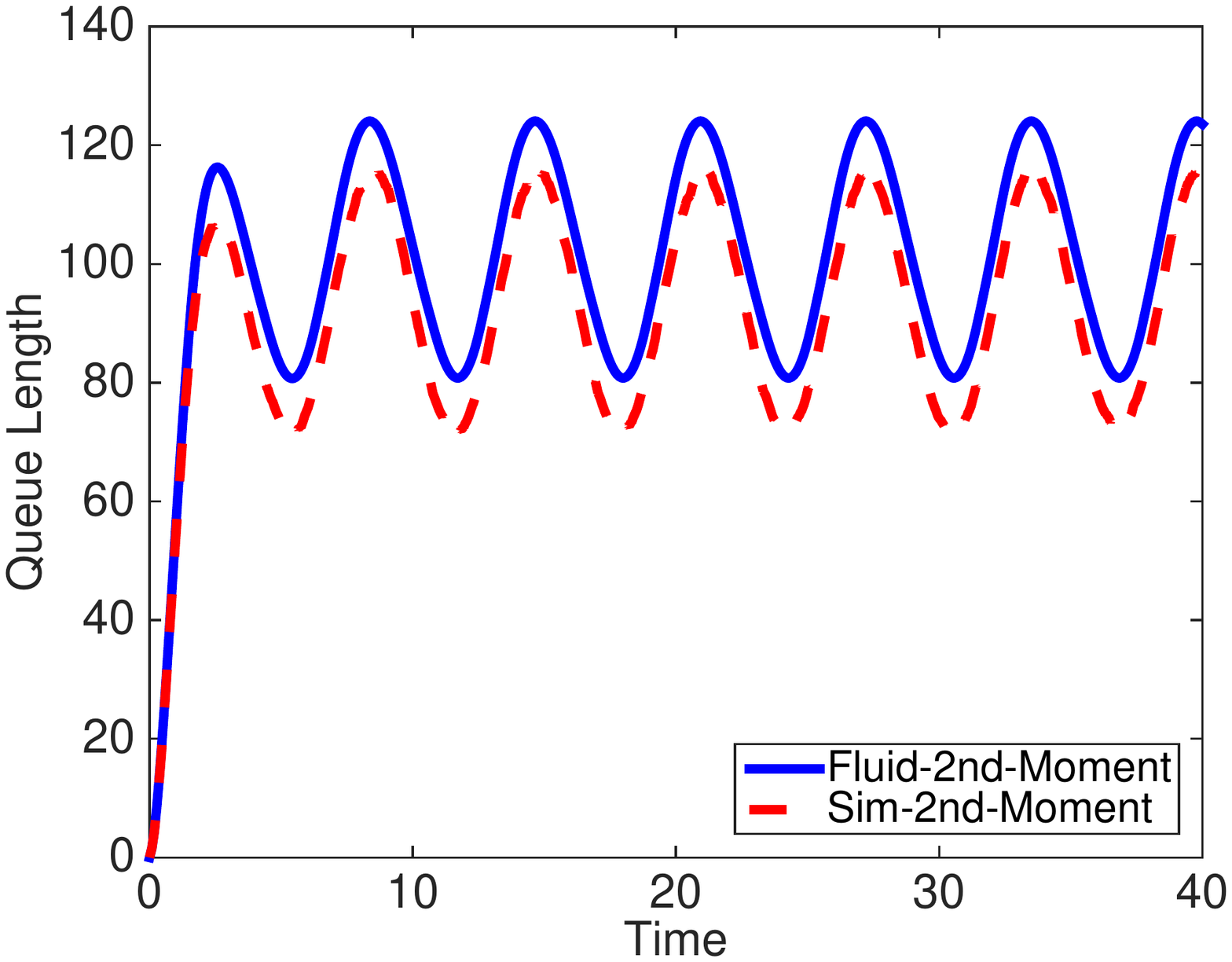}
              \vspace{-.8in}
            \caption{{\small Second Moment }}
            \label{fig:2}
        \end{subfigure}
        \begin{subfigure}[b]{0.475\textwidth}
            \centering   \vspace{-.75in}
            \includegraphics[scale = .3]{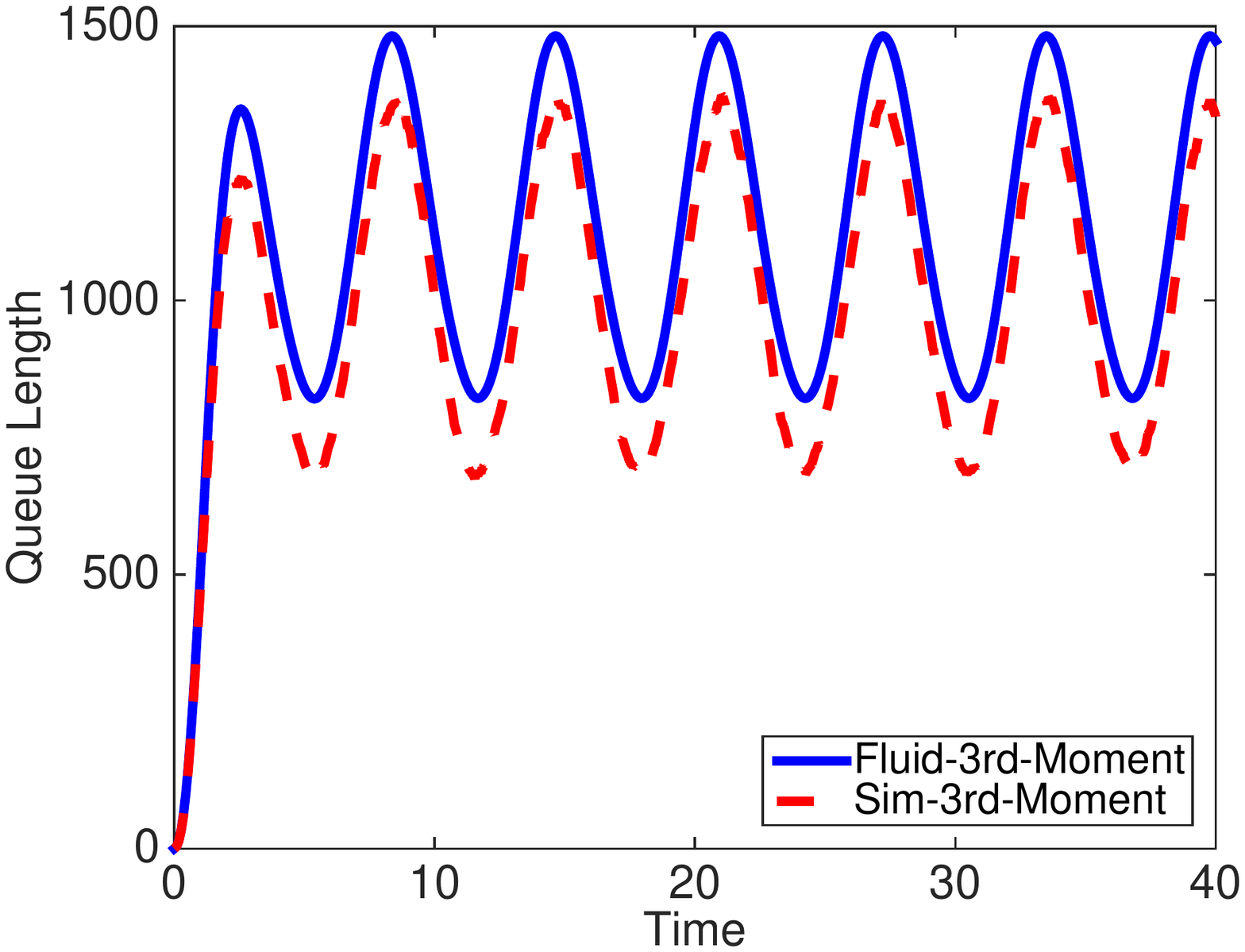}
              \vspace{-.8in}
            \caption{{\small Third Moment }}
            \label{fig:3}
        \end{subfigure}
        \hfill
        \begin{subfigure}[b]{0.475\textwidth}
            \centering   \vspace{-.75in}
            \includegraphics[scale = .3]{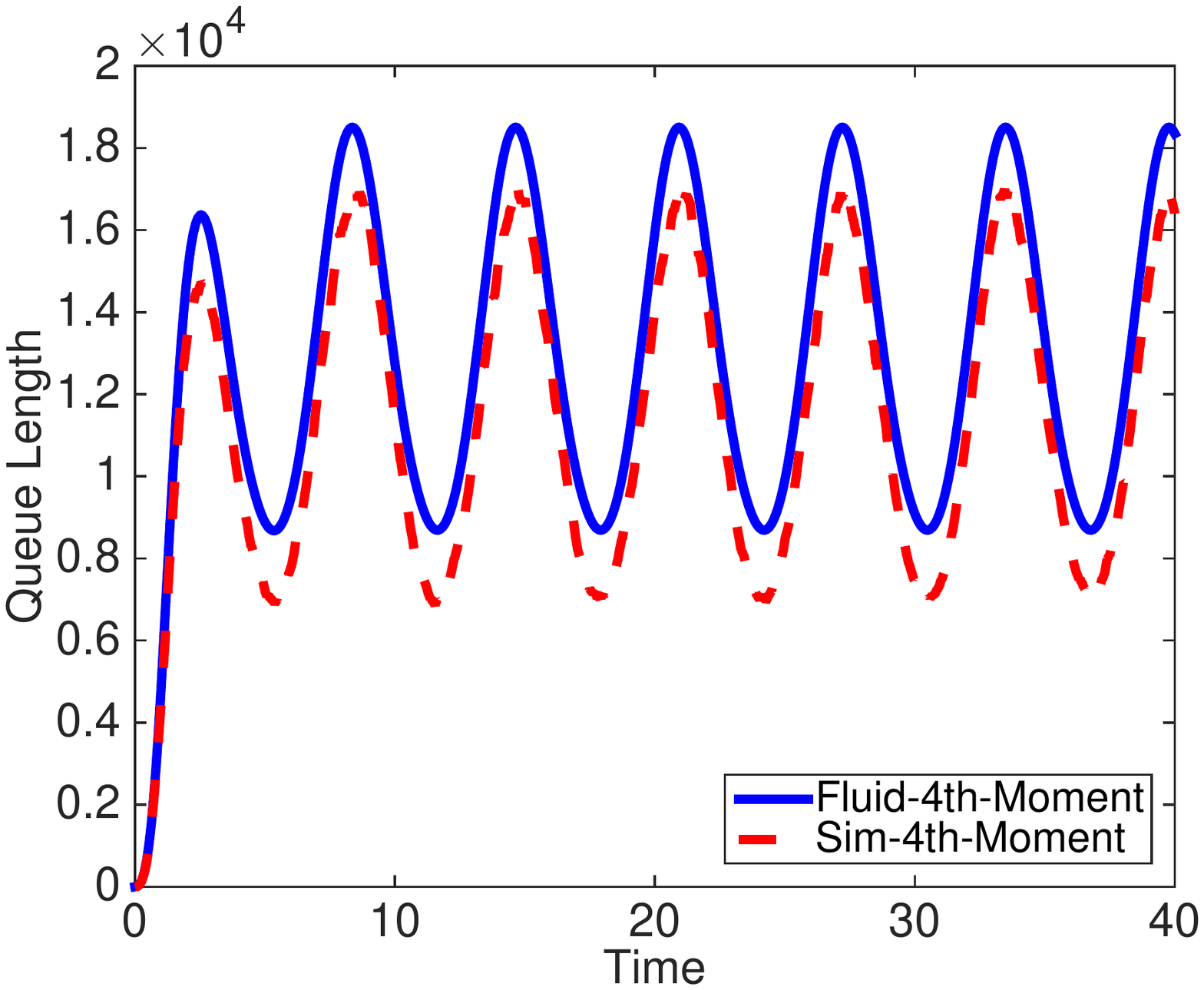}
              \vspace{-.8in}
            \caption{{\small Fourth Moment}}
            \label{fig:4}
        \end{subfigure}
        \caption
        { $\lambda(t) = 10 + 2 \cdot \sin( t) $, $\mu = 1$, $\theta = 2$, $Q(0) = 0$, $c=10$.
 }
        \label{fig:1to9pt2}
    \end{figure}

In this section we describe numerical results for approximating the moments of the Erlang-A queue and examine them relative to our findings.  In Figures \ref{fig:1to9pt1} and \ref{fig:1to9pt2}, we show the first four moments of the Erlang-A queue and their respective fluid approximations for cases of $\theta < \mu$ and $\theta > \mu$, respectively. In these plots, we take the arrival rate at time $t \geq 0$ to be $\lambda(t) = 10 + 2\sin(t)$. We initialize the queue as empty, and we assume that the queueing system has $c = 10$ servers each with exponential service rate $\mu = 1$. We test two different cases for the abandonment rate: $\theta = 0.5$ and $\theta = 2$. In these settings, we observe that when $\theta < \mu$ the fluid approximations are below their corresponding simulated stochastic values and that when $\theta > \mu$ the fluid values are greater than the simulations, and this matches the statements of Theorems \ref{MeanJensen} and \ref{mMoment}.\\

\begin{figure}[H]
        \centering
          \begin{subfigure}[b]{0.475\textwidth}
            \centering   \vspace{-.75in}
           \includegraphics[scale = .3]{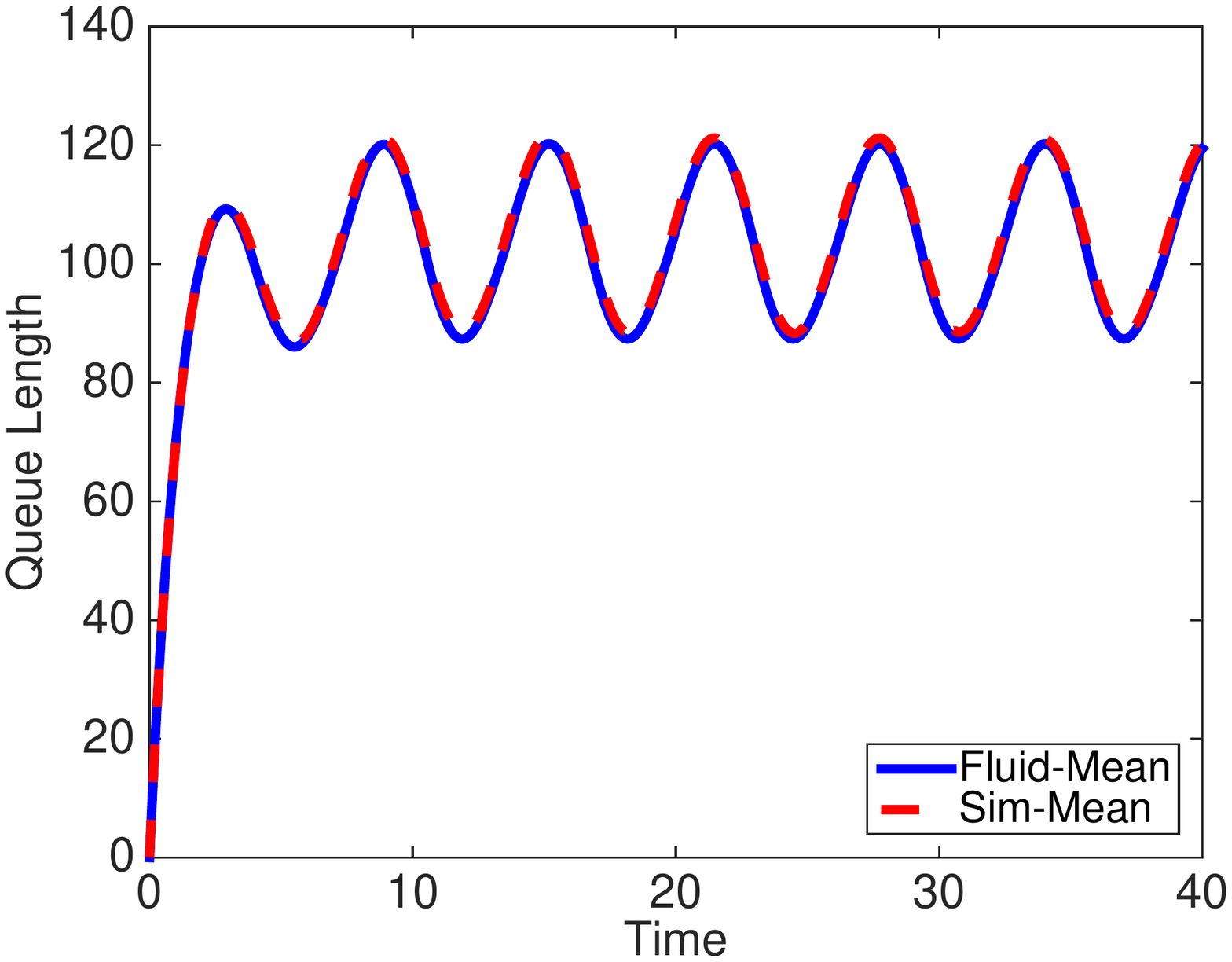}
           \vspace{-.8in}
            \caption{{\small First Moment }}
            \label{fig:1}
        \end{subfigure}
        \hfill
        \begin{subfigure}[b]{0.475\textwidth}
            \centering \vspace{-.75in}
            \includegraphics[scale = .3]{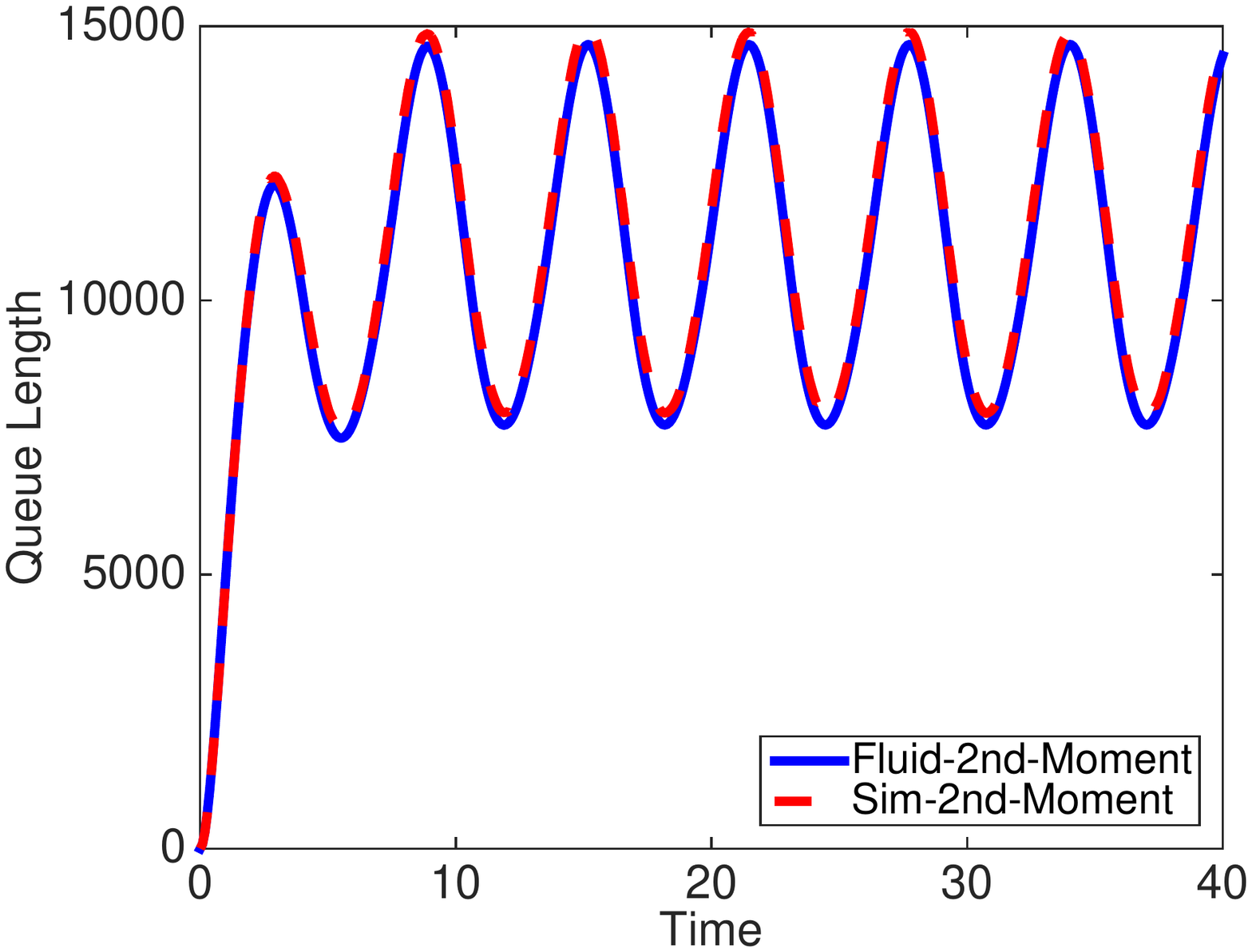}
              \vspace{-.8in}
            \caption{{\small Second Moment }}
            \label{fig:2}
        \end{subfigure}
        \begin{subfigure}[b]{0.475\textwidth}
            \centering   \vspace{-.75in}
            \includegraphics[scale = .3]{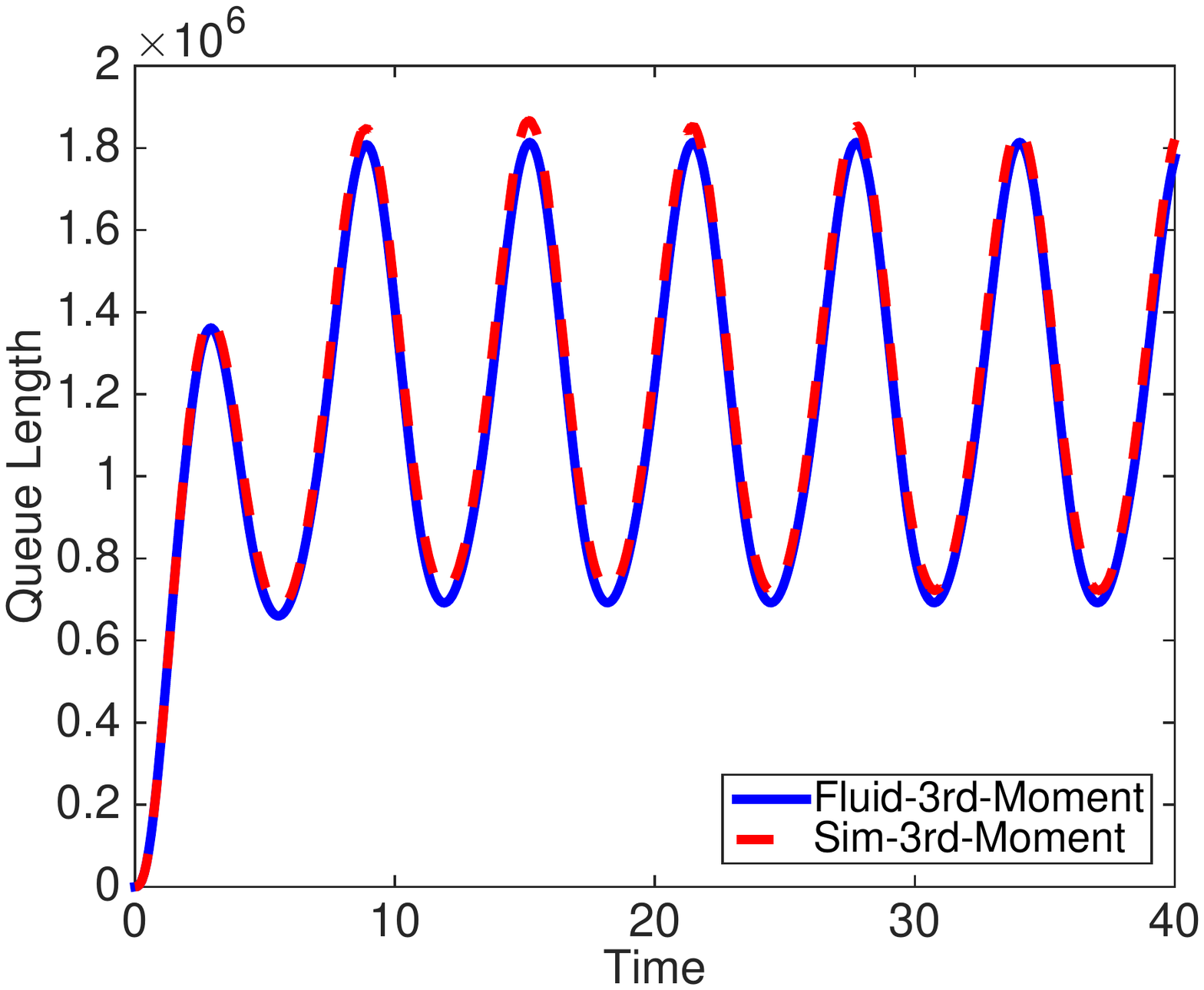}
              \vspace{-.8in}
            \caption{{\small Third Moment }}
            \label{fig:3}
        \end{subfigure}
        \hfill
        \begin{subfigure}[b]{0.475\textwidth}
            \centering   \vspace{-.75in}
            \includegraphics[scale = .3]{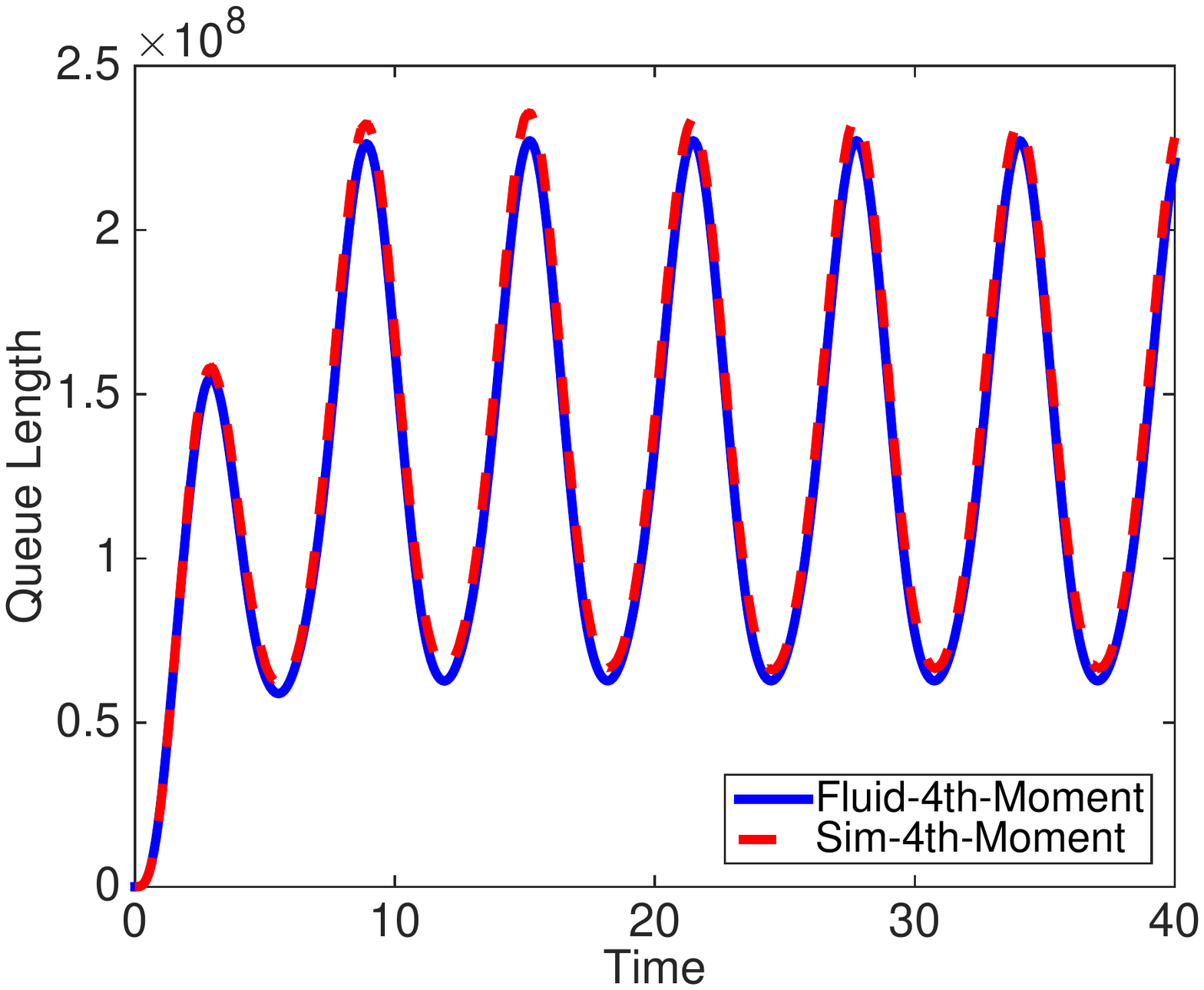}
              \vspace{-.8in}
            \caption{{\small Fourth Moment }}
            \label{fig:4}
        \end{subfigure}
         \caption
        { $\lambda(t) = 100 + 20 \cdot \sin( t) $, $\mu = 1$,  $\theta = 0.5$, $Q(0) = 0$, $c=100$.
 }
        \label{fig:11to19pt1}
     \end{figure}

\begin{figure}[H]
        \begin{subfigure}[b]{0.475\textwidth}
            \centering   \vspace{-.75in}
           \includegraphics[scale = .3]{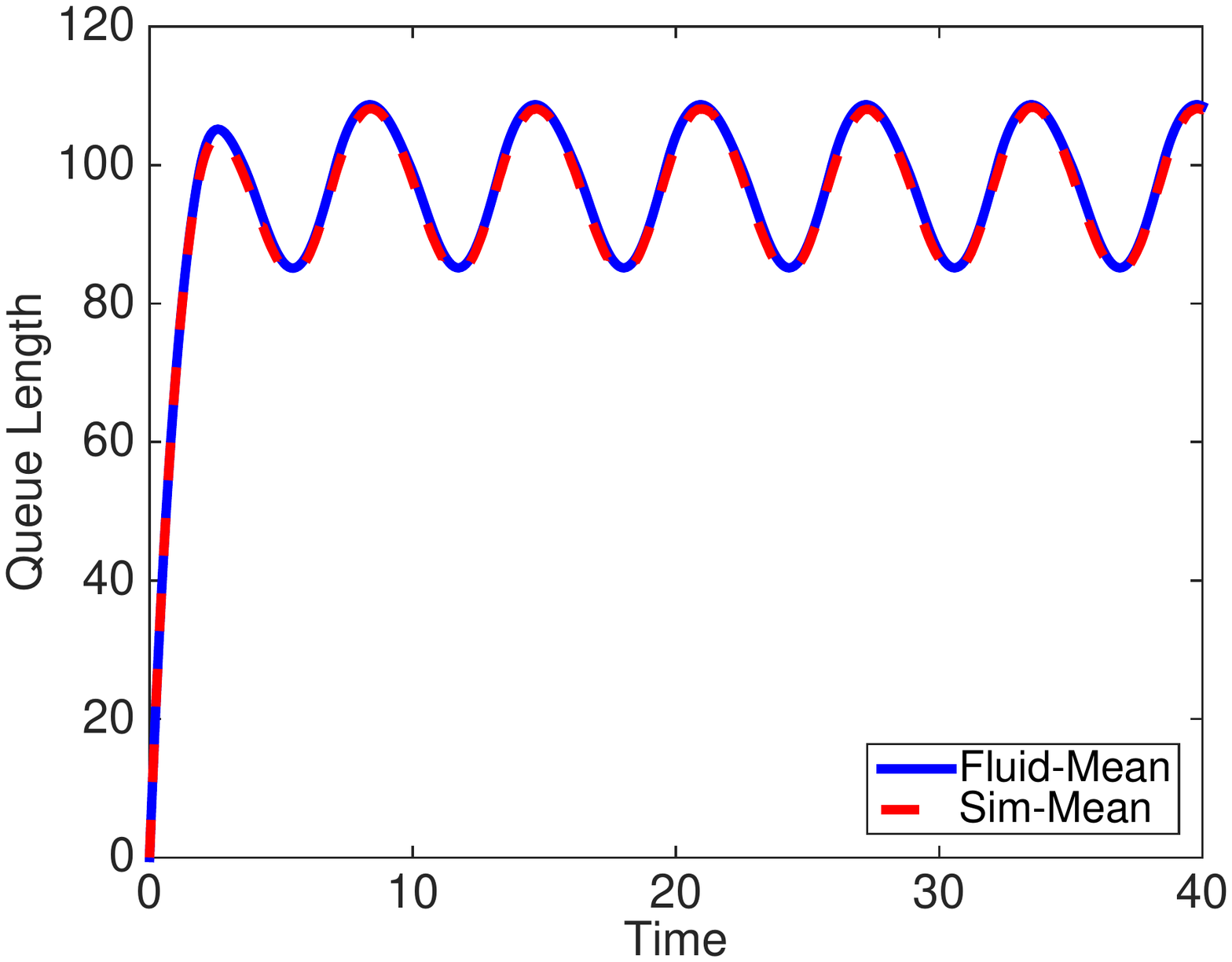}
           \vspace{-.8in}
            \caption{{\small First Moment }}
            \label{fig:1}
        \end{subfigure}
        \hfill
        \begin{subfigure}[b]{0.475\textwidth}
            \centering \vspace{-.75in}
            \includegraphics[scale = .3]{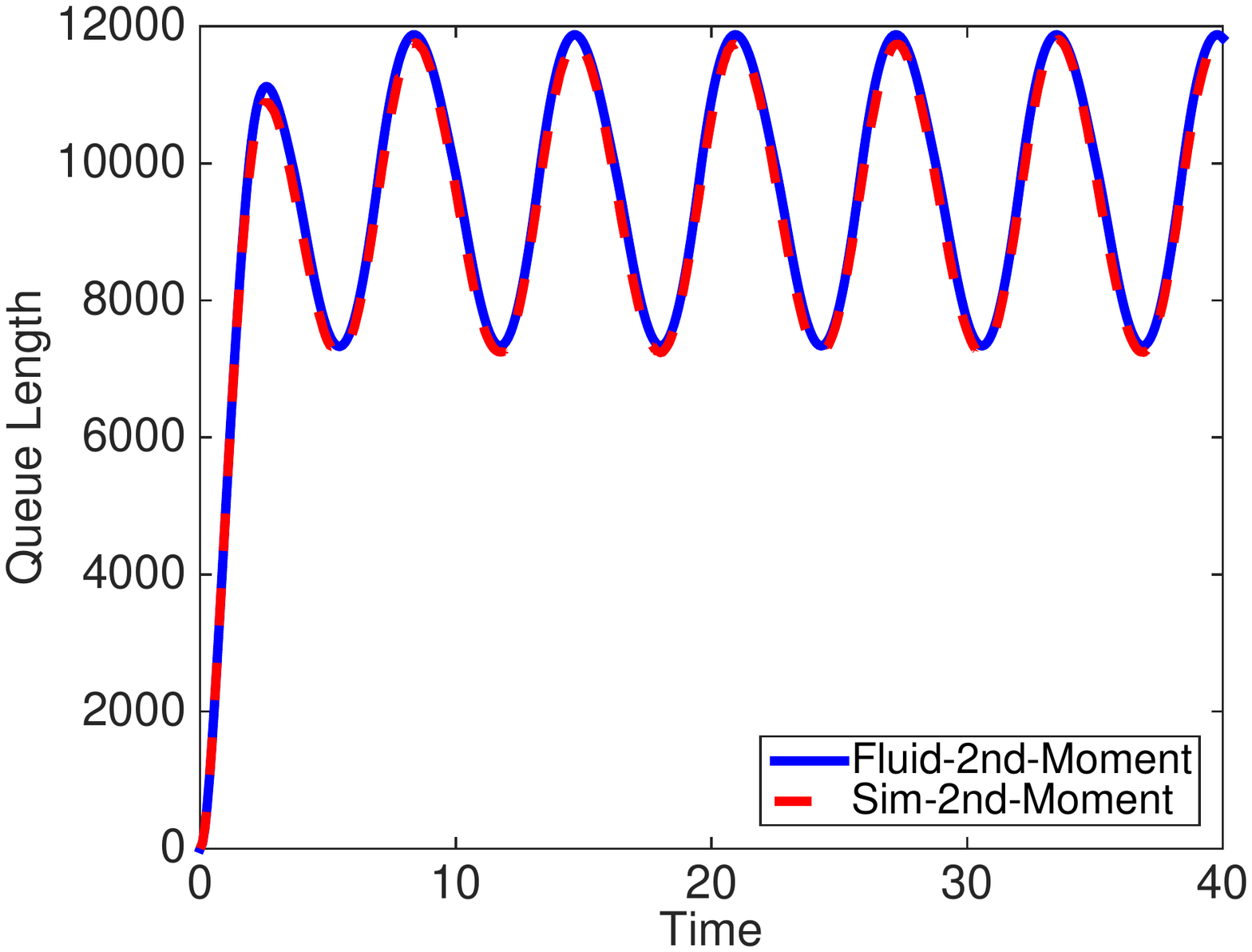}
              \vspace{-.8in}
            \caption{{\small Second Moment }}
            \label{fig:2}
        \end{subfigure}
        \begin{subfigure}[b]{0.475\textwidth}
            \centering   \vspace{-.75in}
            \includegraphics[scale = .3]{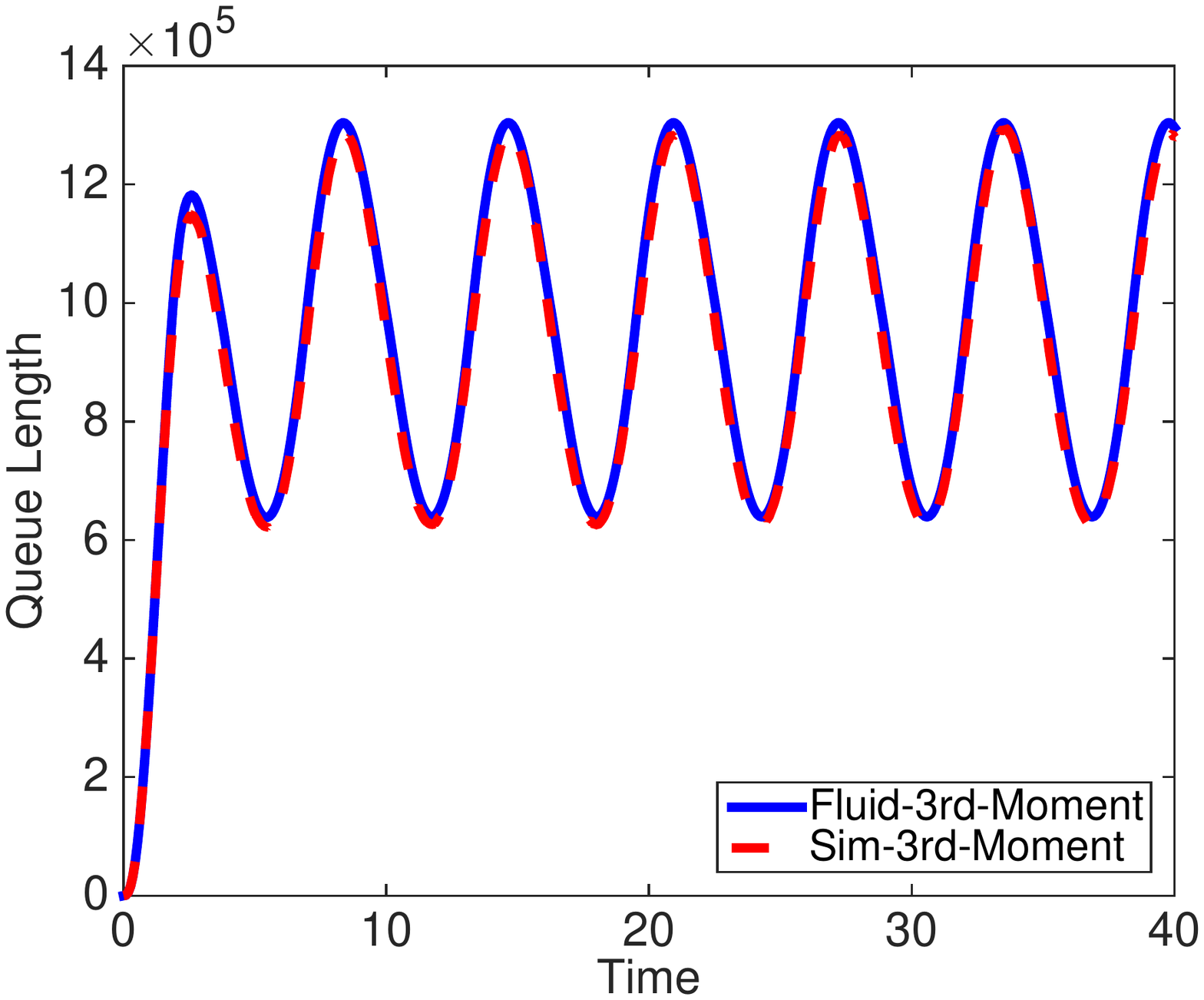}
              \vspace{-.8in}
            \caption{{\small Third Moment }}
            \label{fig:3}
        \end{subfigure}
        \hfill
        \begin{subfigure}[b]{0.475\textwidth}
            \centering   \vspace{-.75in}
            \includegraphics[scale = .3]{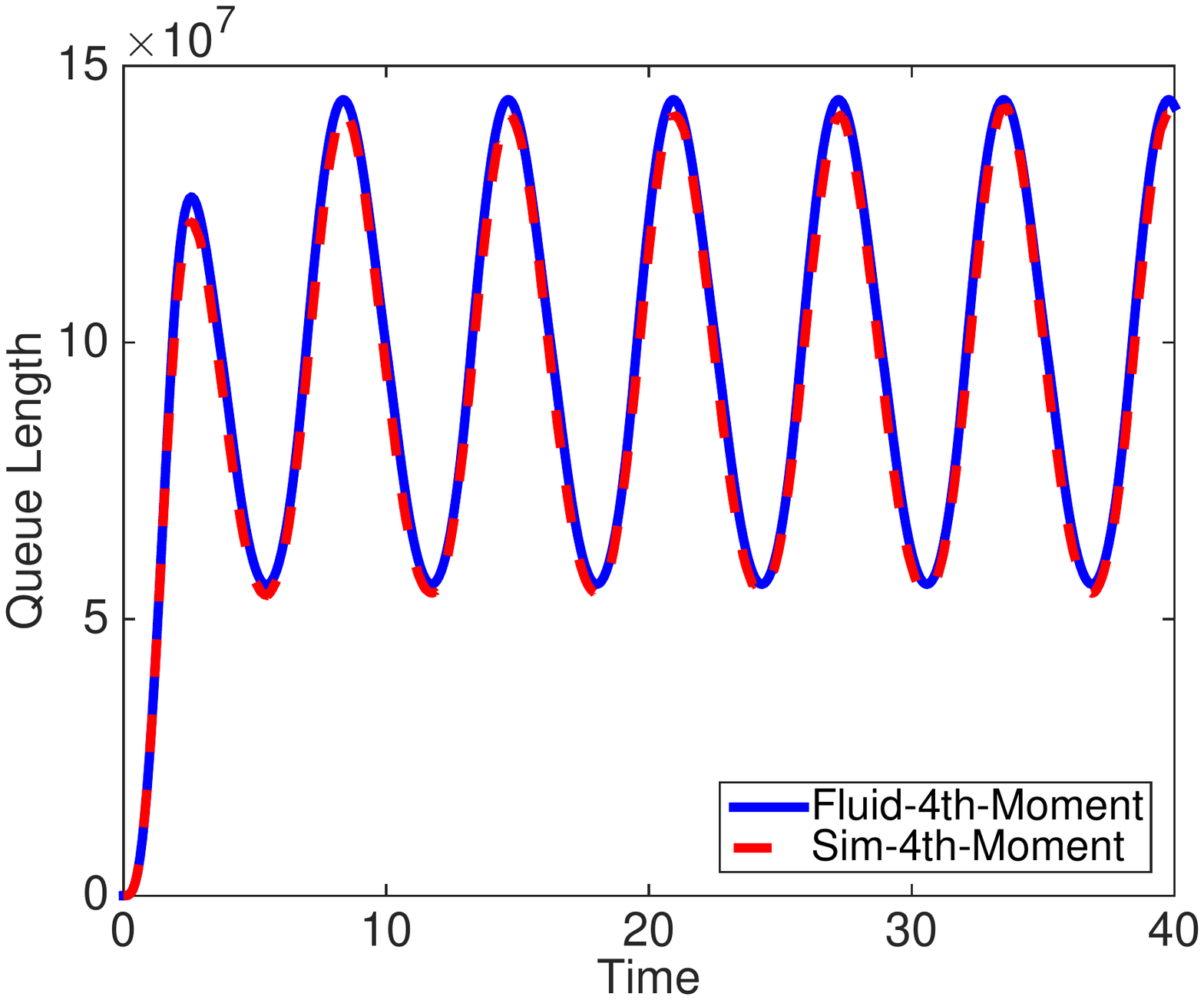}
              \vspace{-.8in}
            \caption{{\small Fourth Moment }}
            \label{fig:4}
        \end{subfigure}
        \caption
        { $\lambda(t) = 100 + 20 \cdot \sin( t) $, $\mu = 1$,  $\theta = 2$, $Q(0) = 0$, $c=100$.
 }
        \label{fig:11to19pt2}
    \end{figure}

We observe the same relationships in Figures~\ref{fig:11to19pt1} and \ref{fig:11to19pt2}. For these plots we instead set $\lambda(t) = 100 + 20\sin(t)$ and $c = 100$ and otherwise use the same values as for  Figures \ref{fig:1to9pt1} and \ref{fig:1to9pt2}. With this increase in the arrival intensity and the number of servers, we see that the gaps between the fluid approximations and the simulations are again present,  albeit proportionally smaller.

\vspace{-.1in}

\section{Inequalities and Characterizations for Generating Functions of the Erlang-A Queue} \label{MGFsec}

Building on what we have found for the moments of the Erlang-A, we can provide similar inequalities for the moment generating function and the cumulant generating function again through convexity in the differential equations for the fluid approximations. We provide these inequalities in Subsections~\ref{MGFJensen} and \ref{CGFJensen}, respectively. In doing so, we find forms for the fluid approximations that we can interpret in terms of expectations of other random quantities. Through these recognitions, we characterize the fluid approximations. We describe these representations for systems in steady-state in Subsection~\ref{steadychar} and for nonstationary systems in Subsection~\ref{nonstatchar}. We conclude this section with  a variety of demonstrations of these results through empirical experiments in Subsection~\ref{mgfnumres}.

\subsection{An Inequality for the Moment Generating Function of the Erlang-A Queue} \label{MGFJensen}

  Using the functional forward equations \citet{massey2013gaussian}, we can show that the moment generating function for the Erlang-A queue satisfies the following partial differential equation
  \begin{eqnarray}\label{mgfpde}
     \updot{\mathrm{E}} \left[e^{\alpha \cdot Q(t)} \right]  &=&  \lambda(t) \cdot (e^{\alpha} -1) \cdot \mathrm{E}\left[e^{\alpha \cdot Q(t)} \right] +  \theta \cdot (e^{-\alpha} -1) \cdot \mathrm{E}\left[ Q(t) \cdot e^{\alpha \cdot Q(t)} \right] \\
     &&- ( \theta - \mu) \cdot (e^{-\alpha} -1) \cdot \mathrm{E}\left[ ( Q(t) \wedge c) \cdot e^{\alpha \cdot Q(t)}\right] \\
     &=&  \lambda(t) \cdot (e^{\alpha} -1) \cdot \mathrm{E}\left[e^{\alpha \cdot Q(t)} \right] +  \theta \cdot (e^{-\alpha} -1) \cdot \frac{\partial M(t,\alpha)}{\partial \alpha} \\
     &&- ( \theta - \mu) \cdot (e^{-\alpha} -1) \cdot \mathrm{E}\left[ ( Q(t) \wedge c) \cdot e^{\alpha \cdot Q(t)}\right] .
  \end{eqnarray}

Just like the non-autonomous differential equation for the mean in Equation \ref{meanforeqn}, we also cannot directly compute the moment generating function since we do not know the distribution of the queue length a priori.  This is also true for numerical purposes.  Unless we can compute the expectation that includes the minimum function it is impossible to know the moment generating function, except in special cases such as the infinite server queue and some cases of the Erlang-B queue.  Thus, it is useful to obtain approximations that are explicit upper or lower bounds for the moment generating function.  By using Jensen's inequality for concave functions, we can approximate the moment generating function with the following partial differential equation
     \begin{eqnarray}
     \updot{\mathrm{E}} \left[e^{\alpha \cdot Q_f(t)} \right]  &=&  \lambda(t) \cdot (e^{\alpha} -1) \cdot \mathrm{E}\left[e^{\alpha \cdot Q_f(t)} \right] +  \theta \cdot (e^{-\alpha} -1) \cdot \frac{\partial M_f(t,\alpha)}{\partial \alpha} \nonumber\\
     &&- ( \theta - \mu) \cdot (e^{-\alpha} -1) \cdot   \left( \mathrm{E}\left[Q_f(t) \cdot e^{\alpha \cdot Q_f(t)}\right] \wedge \mathrm{E}\left[c \cdot e^{\alpha \cdot Q_f(t)} \right]  \right) \\
    \frac{\partial M_f(t,\alpha) }{\partial t} &=&  \lambda(t) \cdot (e^{\alpha} -1) \cdot M_f(t,\alpha) +  \theta \cdot (e^{-\alpha} -1) \cdot \frac{\partial M_f(t,\alpha)}{\partial \alpha} \nonumber\\
     &&- ( \theta - \mu) \cdot (e^{-\alpha} -1) \cdot   \left( \frac{\partial M_f(t,\alpha) }{\partial \alpha} \wedge c \cdot M_f(t,\alpha) \right) \label{fluidMGFdef}.
  \end{eqnarray}

The following theorem determines exactly when $\mathrm{E}\left[e^{\alpha \cdot Q_f(t)} \right]$  is a lower or upper bound for the exact moment generating function of the Erlang-A queue.

\begin{theorem}\label{mgfineqthm}
For the Erlang-A queue, if $Q(0) = Q_{f}(0)$, then $\mathrm{E}\left[e^{\alpha \cdot Q(t)} \right]  \geq  \mathrm{E}\left[e^{\alpha \cdot Q_f(t)} \right] $ when $\theta < \mu$, $\mathrm{E}\left[e^{\alpha \cdot Q(t)} \right] \leq  \mathrm{E}\left[e^{\alpha \cdot Q_f(t)} \right] $ when $\theta > \mu$, and $\mathrm{E}\left[e^{\alpha \cdot Q(t)} \right]  = \mathrm{E}\left[e^{\alpha \cdot Q_f(t)} \right] $ when $\theta = \mu$.

\begin{proof}
 If we take the difference of the two partial differential equations, we obtain the following
     \begin{eqnarray*}
     \updot{\mathrm{E}} \left[e^{\alpha \cdot Q(t)} \right]  -  \updot{\mathrm{E}} \left[e^{\alpha \cdot Q_f(t)} \right]   &=&  \lambda(t) \cdot (e^{\alpha} -1) \cdot \mathrm{E}\left[e^{\alpha \cdot Q(t)} \right] +  \theta \cdot (e^{-\alpha} -1) \cdot \mathrm{E}\left[ Q(t) \cdot e^{\alpha \cdot Q(t)} \right] \\
     &&- ( \theta - \mu) \cdot (e^{-\alpha} -1) \cdot \mathrm{E}\left[ ( Q(t) \wedge c) \cdot e^{\alpha \cdot Q(t)}\right] \\
&&-  \lambda(t) \cdot (e^{\alpha} -1) \cdot \mathrm{E}\left[e^{\alpha \cdot Q_f(t)} \right] -  \theta \cdot (e^{-\alpha} -1) \cdot \mathrm{E}\left[Q_f(t) \cdot e^{\alpha \cdot Q_f(t)}\right]  \\
     &&+ ( \theta - \mu) \cdot (e^{-\alpha} -1) \cdot   \left( \mathrm{E}\left[Q_f(t) \cdot e^{\alpha \cdot Q_f(t)}\right] \wedge \mathrm{E}\left[c \cdot e^{\alpha \cdot Q_f(t)} \right]  \right) \\
     &=&  \lambda(t) \cdot (e^{\alpha} -1) \cdot \left(  \mathrm{E}\left[e^{\alpha \cdot Q(t)} \right] - \mathrm{E}\left[e^{\alpha \cdot Q_f(t)} \right]  \right) \\
     &&+  \theta \cdot (e^{-\alpha} -1) \cdot \left( \mathrm{E}\left[ Q(t) \cdot e^{\alpha \cdot Q(t)} \right] - \mathrm{E}\left[Q_f(t) \cdot e^{\alpha \cdot Q_f(t)}\right] \right) \\
     &&- ( \theta - \mu) \cdot (e^{-\alpha} -1) \cdot \mathrm{E}\left[ ( Q(t) \wedge c) \cdot e^{\alpha \cdot Q(t)}\right]  \\
     &&+ ( \theta - \mu) \cdot (e^{-\alpha} -1) \cdot \left( \mathrm{E}\left[Q_f(t) \cdot e^{\alpha \cdot Q_f(t)}\right] \wedge \mathrm{E}\left[c \cdot e^{\alpha \cdot Q_f(t)} \right]  \right)  .
  \end{eqnarray*}
  Now by exploiting the positive scalability property and the concavity of the minimum function, we have by Jensen's inequality that
  \begin{eqnarray*}
   \mathrm{E}\left[ ( Q(t) \wedge c) \cdot e^{\alpha \cdot Q(t)} \right]   &=&   \mathrm{E}\left[ \left( Q(t) \cdot e^{\alpha \cdot Q(t)}  \wedge c \cdot e^{\alpha \cdot Q(t)} \right)  \right] \\ & \leq & \left( \mathrm{E}\left[Q_f(t) \cdot e^{\alpha \cdot Q_f(t)}\right] \wedge \mathrm{E}\left[c \cdot e^{\alpha \cdot Q_f(t)} \right]  \right) .
  \end{eqnarray*}
  Thus, we have when $\theta < \mu$ that
  \begin{equation}
     \updot{\mathrm{E}} \left[e^{\alpha \cdot Q(t)} \right]  -  \updot{\mathrm{E}} \left[e^{\alpha \cdot Q_f(t)} \right]  \geq 0 ,
     \end{equation}
     when $\theta > \mu$
       \begin{equation}
     \updot{\mathrm{E}} \left[e^{\alpha \cdot Q(t)} \right]  -  \updot{\mathrm{E}} \left[e^{\alpha \cdot Q_f(t)} \right]  \leq 0 ,
     \end{equation}
      and finally when $\theta = \mu$,
          \begin{equation}
     \updot{\mathrm{E}} \left[e^{\alpha \cdot Q(t)} \right]  -  \updot{\mathrm{E}} \left[e^{\alpha \cdot Q_f(t)} \right]  = 0
     \end{equation}
     since they solve the same partial differential equation.  This completes our proof.
\end{proof}
\end{theorem}

As with the moments, we can observe these relationships occurring in numerical experiments. We provide figures demonstrating this in Subsection~\ref{mgfnumres}.

\subsection{An Inequality for the Cumulant Moment Generating Function of the Erlang-A Queue}\label{CGFJensen}

As a consequence of the findings for the moment generating function, we can also provide similar inequalities for the cumulant moment generating function. Using Equation~\ref{mgfpde}, we have
  \begin{eqnarray}
     \updot{\log\left( \mathrm{E} \left[e^{\alpha \cdot Q(t)} \right] \right) }  &\equiv& \frac{\partial}{\partial t} \log\left( \mathrm{E} \left[e^{\alpha \cdot Q(t)} \right] \right)
     =  \frac{\updot{\mathrm{E}} \left[e^{\alpha \cdot Q(t)} \right] }{ \mathrm{E} \left[e^{\alpha \cdot Q(t)} \right] } \\
     &=& \lambda(t) \cdot (e^{\alpha} -1) +  \theta \cdot (e^{-\alpha} -1) \cdot \frac{ \mathrm{E}\left[ Q(t) \cdot e^{\alpha \cdot Q(t)} \right] }{ \mathrm{E} \left[e^{\alpha \cdot Q(t)} \right] } \nonumber \\
     &&- ( \theta - \mu) \cdot (e^{-\alpha} -1) \cdot \frac{ \mathrm{E}\left[ ( Q(t) \wedge c) \cdot e^{\alpha \cdot Q(t)}\right] }{ \mathrm{E} \left[e^{\alpha \cdot Q(t)} \right] }\\
     &=&  \lambda(t) \cdot (e^{\alpha} -1) +  \theta \cdot (e^{-\alpha} -1) \cdot \frac{\partial G(t,\alpha)}{\partial \alpha} \nonumber \\
     &&- ( \theta - \mu) \cdot (e^{-\alpha} -1) \cdot \frac{ \mathrm{E}\left[ ( Q(t) \wedge c) \cdot e^{\alpha \cdot Q(t)}\right] }{\mathrm{E} \left[e^{\alpha \cdot Q(t)} \right]} .
  \end{eqnarray}
Like for the MGF, we note that we cannot compute the cumulant moment generating function directly without knowing the distribution of the queue length. So, by again applying Jensen's inequality, we can describe the fluid approximation as follows.
     \begin{eqnarray}
    \updot{\log\left( \mathrm{E} \left[e^{\alpha \cdot Q_f(t)} \right] \right) }   &=&  \lambda(t) \cdot (e^{\alpha} -1) +  \theta \cdot (e^{-\alpha} -1) \cdot \frac{\partial G_f(t,\alpha)}{\partial \alpha} \nonumber\\
     &&- ( \theta - \mu) \cdot (e^{-\alpha} -1) \cdot   \left( \frac{ \mathrm{E}\left[Q_f(t) \cdot e^{\alpha \cdot Q_f(t)}\right] \wedge \mathrm{E}\left[c \cdot e^{\alpha \cdot Q_f(t)} \right] }{ \mathrm{E} \left[e^{\alpha \cdot Q(t)} \right] }  \right) \quad\\
    \frac{\partial G_f(t,\alpha) }{\partial t} &=&  \lambda(t) \cdot (e^{\alpha} -1) +  \theta \cdot (e^{-\alpha} -1) \cdot \frac{\partial G_f(t,\alpha)}{\partial \alpha} \nonumber\\
     &&- ( \theta - \mu) \cdot (e^{-\alpha} -1) \cdot   \left( \frac{\partial G(t,\alpha) }{\partial \alpha} \wedge c  \right) .\label{cumulantMGF}
  \end{eqnarray}

Using this observation and our approach in finding the inequalities for the moment generating function, we find the equivalent inequalities for the cumulant moment generating function in the following corollary.

\begin{corollary}
For the Erlang-A queue, if $Q(0) = Q_{f}(0)$, then $\log\left(\mathrm{E}\left[e^{\alpha \cdot Q(t)} \right]\right)  \geq  \log\left(\mathrm{E}\left[e^{\alpha \cdot Q_f(t)} \right]\right) $ when $\theta < \mu$, $\log\left(\mathrm{E}\left[e^{\alpha \cdot Q(t)} \right]\right) \leq  \log\left(\mathrm{E}\left[e^{\alpha \cdot Q_f(t)} \right] \right)$ when $\theta > \mu$, and $\log\left(\mathrm{E}\left[e^{\alpha \cdot Q(t)} \right]\right)  = \log\left(\mathrm{E}\left[e^{\alpha \cdot Q_f(t)} \right] \right)$ when $\theta = \mu$.

\begin{proof}
The proof follows from the same argument that was given in Theorem~\ref{mgfineqthm} and the fact that the log function is strictly increasing.
\end{proof}
\end{corollary}

\subsection{Characterization of the Moment Generating Function in Steady-State}\label{steadychar}

From what we have observed for the moment generating function, we can derive an exact representation for the fluid approximation of the moment generating function in steady-state. We assume a stationary arrival rate $\lambda > 0$. We will  investigate the stationary fluid approximation differential equations in a casewise manner based on the relationship of $\lambda$ and the system's service parameters. To do so, we begin with a lemma bounding the fluid approximation of the mean.

\begin{lemma}\label{cLemma}
Suppose that $\lambda$ is constant. If $\lambda < c\mu $, then $E[Q_f(\infty)] < c $.  Moreover, if $\lambda \geq c\mu $, then $E[Q_f(\infty)] \geq c $.
\begin{proof}
We will prove this by contradiction.  For the first part, we assume that $E[Q_f(\infty)] \geq c $.  Now by using the differential equation for the mean in steady state, we have that
\begin{eqnarray*}
0 &=& \lambda - \mu \cdot (E[Q_f(\infty)] \wedge c ) - \theta \cdot ( E[Q_f(\infty)] - c)^+ \\
&=& \lambda - \mu \cdot c  - \theta ( E[Q_f(\infty)] - c)^+ .
\end{eqnarray*}
Since we assumed that $E[Q_f(\infty)] \geq c $, then this yields the following inequality
\begin{eqnarray*}
 \lambda &\geq & c\mu ,
\end{eqnarray*}
which yields a contradiction.
For the second case, where we assume that $\lambda \geq c \mu$ and $E[Q_f(\infty)] < c $, then by the same differential equation we have that
\begin{eqnarray*}
 \lambda &=&  \mu \cdot (E[Q_f(\infty)] \wedge c ) + \theta \cdot ( E[Q_f(\infty)] - c)^+ \\
 &=& \mu \cdot (E[Q_f(\infty)] \wedge c ) \\
 &=& c \mu + \mu \cdot (E[Q_f(\infty)] - c ) \\
 &< & c \mu ,
\end{eqnarray*}
which yields another contradiction.
\end{proof}
\end{lemma}

We now begin characterizing the fluid approximations with our first case, $\lambda \geq c\mu$, in the following proposition.

\begin{proposition}\label{fluidSteadyMGF}
If $\lambda \geq c\mu$, then in steady-state we have that
   \begin{eqnarray}
   \frac{\partial M_f(\infty,\alpha)}{\partial \alpha} &=&  \frac{\lambda \cdot (e^{\alpha} -1) +  ( \theta - \mu) \cdot ( 1 - e^{-\alpha}) \cdot  c }{\theta \cdot ( 1 - e^{-\alpha}) } \cdot  M_f(\infty,\alpha)
  \end{eqnarray}
which yields a solution of
     \begin{eqnarray}
   M_f(\infty,\alpha)   &=& e^{\frac{\alpha \cdot (\theta - \mu) \cdot c + \lambda \cdot (e^\alpha -1)}{\theta }}
  \end{eqnarray}
  for $\alpha \in \mathbb{R}$.
  \begin{proof}
  To find the partial differential equation, we use functional cumulant bound for any non-decreasing function $h(\cdot)$ (which can be seen as a form of the FKG inequality),
  \begin{equation}
  \frac{\mathrm{E}[ h(X) \cdot e^{\alpha \cdot X}] }{\mathrm{E}[ e^{\alpha \cdot X}] } \geq \mathrm{E}[ h(X) ] .
  \end{equation}
  In the case that $\lambda \geq c\mu$ we have that $\E{Q_f(t)} \geq c$ in steady-state by Lemma~\ref{cLemma}, and so we know how to evaluate the minimum in the fluid equation. Thus, we have that the derivative of $G_f(\alpha) = \log(M_f(\infty, \alpha))$ with respect to $\alpha$ is
\begin{align}\label{logMGFode}
\frac{\mathrm{d}G_f(\alpha)}{\mathrm{d}\alpha}
=
\frac{\lambda(e^{\alpha} - 1 ) + c(\theta - \mu)(1 - e^{-\alpha})}{\theta(1 - e^{-\alpha})} = \frac{\lambda e^\alpha}{\theta} + \frac{c(\theta - \mu)}{\theta}
\end{align}
where here we have used the identity $e^x = \frac{e^x - 1}{1 - e^{-x}}$, which can be observed by multiplying each side of the equation by $1 - e^{-x}$. Because the MGF is equal to 1 when $\alpha = 0$, we also have that $G_f(0) = 0$. Using this initial condition and integrating left and right sides of Equation~\ref{logMGFode} with respect to $\alpha$, we find that
$$
G_f(\alpha) = \frac{\lambda(e^\alpha - 1) + c\alpha(\theta - \mu)}{\theta}
$$
and since $M_f(\infty, \alpha) = e^{G_f(\alpha)}$, we attain the stated result.
  \end{proof}
\end{proposition}

We can now observe that the fluid approximation is equivalent in distribution to a Poisson random variable shifted by $\gamma \equiv \frac{c(\theta - \mu)}{\theta}$, as the moment generation function for the Poisson distribution is $e^{\beta (e^{\alpha} - 1)}$, where $\beta$ is the rate of arrival and $\alpha$ is the space parameter of the MGF. This gives rise to the following.

\begin{theorem}\label{sandwichthm}
For the Erlang-A queue with $\lambda \geq c\mu$ and $m \in \mathbb{Z}^+$, if $\theta > \mu$
$$
\E{(Q_f(\infty) - \gamma)^m} \leq \E{(Q(\infty))^m} \leq \E{(Q_f(\infty))^m}
$$
and if $\theta < \mu$
$$
\E{(Q_f(\infty))^m} \leq \E{(Q(\infty))^m} \leq \E{(Q_f(\infty) - \gamma)^m}
$$
where $\gamma = \frac{c(\theta - \mu)}{\theta}$.
\begin{proof}
From Proposition~\ref{fluidSteadyMGF}, we have that the fluid approximation of the MGF in steady-state is
$$
M_f(\infty, \alpha) = e^{\frac{\lambda(e^\alpha - 1) + c\alpha(\theta - \mu)}{\theta}} = \E{e^{\alpha(\Gamma + \gamma)}}
$$
where $\Gamma \sim \mathrm{Pois}\left(\frac{\lambda}{\theta}\right)$ and $\gamma = \frac{c(\theta - \mu)}{\theta}$. From the uniqueness of MGF's, we have that
$$
\E{(Q_f(\infty))^m} = \E{(\Gamma + \gamma)^m}
$$
for all $m \in \mathbb{Z}^+$. Now, recall that for an $M/M/\infty$ queue with arrival rate $\lambda$ and service rate $\theta$, the stationary distribution is that of a Poisson random variable with rate parameter $\frac{\lambda}{\theta}$. So, we can think of $\Gamma$ as representing the steady-state distribution of an infinite server queue with Poisson arrival rate $\lambda$ and exponential service rate $\theta$.\\

Suppose now that $\theta > \mu$. Then, by Theorem~\ref{mMoment} and our preceding observation, we have that $\E{(Q(\infty))^m} \leq \E{(\Gamma + \gamma)^m} $. Additionally, by comparing the steady-state infinite server queue representation of $\Gamma$ to $Q(\infty)$, we can further observe that $\E{(Q(\infty))^m} \geq \E{\Gamma^m}$, as for any state $j$ the service rate in $Q(\infty)$ is no more than the service rate in the same state in the $\Gamma$ queueing system. Thus we have that
$$
\E{(Q_f(\infty) - \gamma)^m} = \E{\Gamma^m} \leq \E{(Q(\infty))^m} \leq \E{(\Gamma + \gamma)^m} = \E{(Q_f(\infty))^m}
$$
for all $m \in \mathbb{Z}^+$ whenever $\theta > \mu$.
By symmetric arguments, we can also find that if $\mu > \theta$ then
$$
\E{(Q_f(\infty))^m} = \E{(\Gamma + \gamma)^m} \leq \E{(Q(\infty))^m} \leq \E{\Gamma^m} = \E{(Q_f(\infty) - \gamma)^m}
$$
for all $m \in \mathbb{Z}^+$, as in this case $\gamma = \frac{c(\theta - \mu)}{\theta} < 0$.
\end{proof}
\end{theorem}

\begin{rem}
Note that in Theorem~\ref{mMoment}, we require that $Q(0) = Q_f(0)$ but in this case we have not assumed such a condition. This is because the inequalities in Theorem~\ref{mMoment} hold for all time, and we simply need the relationship to hold in steady-state, which can be seen to occur regardless of initial conditions.
\end{rem}

By knowing the fluid form of moment generating function explicitly as a Poisson distribution, we can also provide exact expressions for the fluid moments and the fluid cumulant moments. These are given in the two following corollaries.

\begin{corollary}
If $\lambda \geq c\mu$, then in steady-state we have that the first $n$ moments have the following steady-state expressions:
  \begin{eqnarray}
\mathrm{E} [Q_{f}^n(\infty)]  &=&  \sum^{n}_{j=0} { n \choose j} \cdot  \left( \frac{c (\theta - \mu)}{\theta} \right)^j \cdot \mathcal{P}_{n-j} \left( \frac{\lambda}{\theta} \right)
\end{eqnarray}
where $\mathcal{P}_m \left( \frac{\lambda}{\theta} \right) $ is the $m^{th}$ Touchard polynomial with parameter $ \frac{\lambda}{\theta} $.

  \begin{proof}
This can be seen by direct use of the Poisson form of the fluid MGF. Let $\Gamma \sim \text{Pois}\left(\frac{\lambda}{\theta}\right)$ and let $\gamma = \frac{c(\theta - \mu)}{\theta}$. Then,
  \begin{eqnarray*}
\mathrm{E} [Q_{f}^n(\infty)]  &=& \mathrm{E} [ (\Gamma + \gamma)^n] \\
 &=&  \sum^{n}_{j=0} { n \choose j } \cdot  \gamma^j \cdot \mathrm{E} \left[ \Gamma^{n-j} \right] \\
 &=& \sum^{n}_{j=0} { n \choose j} \cdot  \gamma^j \cdot \mathcal{P}_{n-j} \left( \frac{\lambda}{\theta} \right) \\
 &=& \sum^{n}_{j=0} { n \choose j} \cdot  \left( \frac{c (\theta - \mu)}{\theta} \right)^j \cdot \mathcal{P}_{n-j} \left( \frac{\lambda}{\theta} \right) .
\end{eqnarray*}
      \end{proof}

\end{corollary}

\begin{corollary}
If $\lambda \geq c\mu$, then in steady-state we have that

\begin{equation}
\frac{\mathrm{d}G_f(\infty,\alpha)}{\mathrm{d}\alpha}
\Bigr|_{\alpha = 0}  = \frac{\lambda }{\theta} + \frac{c(\theta - \mu)}{\theta} = \mathrm{E}[Q_f(\infty)]
\end{equation}
and for $n \in \mathbb{Z}^+$
\begin{equation}
\frac{\mathrm{d}^n G_f(\infty,\alpha)}{\mathrm{d}^n \alpha}
\Bigr|_{\alpha = 0}  = \frac{\lambda}{\theta} = \mathrm{C}^{(n)}[Q_f(\infty)]
\end{equation}
where $\mathrm{C}^{(n)}[Q_f(\infty)]$ is defined as the $n^{th}$ cumulant moment of $Q_f(\infty)$.
\end{corollary}

We now consider the second case, which is $\lambda < c\mu e^{-\alpha}$. Note that this now also requires a relationship involving the space parameter of the moment generating function, $\alpha$. This is less general than the first case, but it allows us to derive Lemma~\ref{ifflemma}.

\begin{lemma}\label{ifflemma}
For $\alpha \geq 0$,
$$
\frac{\partial M_f(\infty, \alpha)}{\partial \alpha} < c M_f(\infty, \alpha)
$$
if and only if $\lambda < c\mu e^{-\alpha}$.
\begin{proof}
To begin, suppose that
$
\frac{\partial M_f(\infty, \alpha)}{\partial \alpha} <  c M_f(\infty, \alpha).
$
Using this information in conjunction with the steady-state form of the partial differential equation for the fluid MGF given in Equation~\ref{fluidMGFdef}, we have that
$$
0
=
\lambda (e^{\alpha} -1) M_f(\infty,\alpha) +  \theta (e^{-\alpha} -1) \frac{\partial M_f(\infty,\alpha)}{\partial \alpha}
     - ( \theta - \mu) (e^{-\alpha} -1) \frac{\partial M_f(\infty, \alpha)}{\partial \alpha}
$$
which simplifies to
$$
\frac{\partial M_f(\infty,\alpha)}{\partial \alpha}
=
\frac{\lambda}{\mu}e^{\alpha} M_f(\infty, \alpha).
$$
Using our assumption, we see that
$$
\frac{\lambda}{\mu}e^{\alpha} M_f(\infty, \alpha) < c M_f(\infty, \alpha)
$$
and this yields that $\lambda < c \mu e^{-\alpha}$, which shows one direction.\\
\\
\indent We now move to showing the opposite direction and instead  assume that $\frac{\partial M_f(\infty, \alpha)}{\partial \alpha} \geq  c M_f(\infty, \alpha)$. In this case, Equation~\ref{fluidMGFdef} is given by
$$
0
=
\lambda (e^{\alpha} -1) M_f(\infty,\alpha) +  \theta (e^{-\alpha} -1) \frac{\partial M_f(\infty,\alpha)}{\partial \alpha}
     - c( \theta - \mu) (e^{-\alpha} -1) M_f(\infty,\alpha)
$$
and this simplifies to
$$
\frac{\partial M_f(\infty, \alpha)}{\partial \alpha}
=
\frac{\lambda(e^{\alpha} - 1 ) + c(\theta - \mu)(1 - e^{-\alpha})}{\theta(1 - e^{-\alpha})}
M_f(\infty, \alpha)
=
\frac{\lambda e^{\alpha}+c(\theta - \mu)}{\theta}
M_f(\infty, \alpha).
$$
Again by use of this case's assumption, we have
$$
\frac{\lambda e^{\alpha}+c(\theta - \mu)}{\theta}
M_f(\infty, \alpha)
\geq
c M_f(\infty, \alpha)
$$
and this now yields
$$
\lambda \geq e^{-\alpha} \left(c \theta - c(\theta - \mu)\right) = c \mu e^{-\alpha},
$$
thus completing the proof.
\end{proof}
\end{lemma}

We can now use this lemma to find an explicit form for the fluid approximation of the steady-state moment generating function when $\lambda < c\mu e^{-\alpha}$.

\begin{proposition}\label{fluidSteadyMGF2}
For $\alpha \geq 0 $, if $\lambda < c\mu e^{-\alpha}$, then in steady-state we have that
     \begin{eqnarray} \label{sslessthan}
   \frac{\partial M_f(\infty,\alpha)}{\partial \alpha} &=&    \frac{\lambda \cdot e^{\alpha}  }{\mu } \cdot  M_f(\infty,\alpha)
  \end{eqnarray}
which yields a solution of
     \begin{eqnarray}
   M_f(\infty,\alpha)   &=& e^{\frac{ \lambda \cdot (e^\alpha -1)}{\mu }} \label{lambdamumgf}
  \end{eqnarray}
  for $\alpha \in \mathbb{R}$.
  \begin{proof}
By Lemma~\ref{ifflemma} and our assumption that $\lambda < c\mu e^{-\alpha}$, we know that $\frac{\partial M_f(\infty, \alpha)}{\partial \alpha} < c M_f(\infty, \alpha)$.  Thus, by observing this in the steady-state MGF equation, we easily obtain the result in Equation \ref{sslessthan}.  Moreover, the solution to Equation \ref{sslessthan} can be easily seen by inserting our proposed solution in and noting that it satisfies our differential equation.  Moreover, the solution is unique by the properties of linear ordinary differential equation theory.
  \end{proof}
\end{proposition}

\begin{rem}
We now pause to note that the $\lambda \geq c\mu e^{-\alpha}$ case of Lemma~\ref{ifflemma} implies Proposition~\ref{fluidSteadyMGF} (and its following consequences) with a weaker assumption. However, because the condition $\lambda \geq c\mu$ does not depend on the choice of $\alpha$ it is more general, and thus we leave those results as stated with that assumption instead of $\lambda \geq c\mu e^{-\alpha}$.
\end{rem}

Here we observe that Equation~\ref{lambdamumgf} is equivalent to the moment generating function of a Poisson random variable with parameter $\frac{\lambda}{\mu}$. Now, by recalling again that the steady-state distribution of a $M/M/\infty$ queue is a Poisson distribution with parameter equal to the arrival rate over the service rate, we find the following inequalities.

\begin{theorem}\label{sandwichthm2}
Let $\lambda < c\mu$ and $m \in \mathbb{Z}^+$. Then, if $\theta > \mu$
\begin{align}
\E{\Gamma_{\theta}^m} \leq \E{Q(\infty)^m} \leq \E{\Gamma_{\mu}^m},
\intertext{and if $\mu > \theta$}
\E{\Gamma_{\mu}^m} \leq \E{Q(\infty)^m} \leq \E{\Gamma_{\theta}^m}
\end{align}
where $\Gamma_{x} \sim \mathrm{Pois}\left(\frac{\lambda}{x}\right)$ for $x > 0$.

\begin{proof}
In each case, the inequality involving $\Gamma_{\mu} \sim \mathrm{Pois}\left(\frac{\lambda}{mu}\right)$ follows directly from Proposition~\ref{fluidSteadyMGF2} and Theorem~\ref{mMoment} via the observation that the fluid form of the moment generating function is equivalent in distribution to that of $\Gamma_\mu$. Here we are using Proposition~\ref{fluidSteadyMGF2} with $\alpha = 0$, and by continuity we know this holds for some ball around 0. This validates the use of the derivatives of the steady-state MGF with respect to $\alpha$ evaluated at $\alpha = 0$ in finding the moments for the fluid approximation. Thus, we are left to prove the inequalities for $\Gamma_\theta \sim \mathrm{Pois}\left(\frac{\lambda}{\theta}\right)$.\\

To do so, let's first note that the stationary distribution of a $M/M/\infty$ queue with service rate $\theta$ is equivalent to that of $\Gamma_\theta$. Suppose now that $\theta > \mu$. Then, any state of such a $M/M/\infty$ queue has a larger rate of departure than the same state in the Erlang-A system. Thus, we have that
$$
\E{\Gamma_{\theta}^m} \leq \E{Q(\infty)^m} \leq \E{\Gamma_{\mu}^m}
$$
for all $m \in \mathbb{Z}^+$. By symmetric arguments in the $\theta < \mu$ case, we complete the proof.
\end{proof}
\end{theorem}

As we did for the case when $\lambda \geq c\mu$, we can use these findings to give explicit expressions for the fluid approximations of the moments and the cumulant moments.

\begin{corollary}
If $\lambda < c\mu$, then in steady-state we have that

\begin{equation}
\frac{\mathrm{d}G_f(\infty, \alpha)}{\mathrm{d}\alpha}
\Bigr|_{\alpha = 0}  = \frac{\lambda }{\mu} = \mathrm{E}[Q_f(\infty)]
\end{equation}
and for $n \in \mathbb{Z}^+$,
\begin{align}
\frac{\mathrm{d}^n G_f(\infty,\alpha)}{\mathrm{d}^n \alpha}
\Bigr|_{\alpha = 0}  &= \frac{\lambda}{\mu} = \mathrm{C}^{(n)}[Q_f(\infty)]\\
\frac{\mathrm{d}^n M_f(\infty,\alpha)}{\mathrm{d}^n \alpha}
\Bigr|_{\alpha = 0}  &= \mathcal{P}_{n} \left( \frac{\lambda}{\mu} \right) = \E{Q_f(\infty)^n}
\end{align}
where $\mathrm{C}^{(n)}[Q_f(\infty)]$ is defined as the $n^{th}$ cumulant moment of $Q_f(\infty)$ and  $\mathcal{P}_m \left( \frac{\lambda}{\mu} \right) $ is the $m^{th}$ Touchard polynomial with parameter $ \frac{\lambda}{\mu} $.
\end{corollary}

\subsection{Characterization of the Nonstationary Moment Generating Function}\label{nonstatchar}

Many scenarios that feature customer abandonments may also feature an arrival process that is nonstationary. To incorporate this, we now incorporate a point process that can be used to approximate any periodic mean arrival pattern, as discussed in \citet{eick1993mt}. Specifically, we define $\lambda(t)$ by a Fourier series: let $\lambda_0$ and $\{(a_k, b_k), k \in \mathbb{Z}^+\}$ be such that
\begin{align}\label{lambdaNonStatDef}
\lambda(t) = \lambda_0 + \sum_{k=1}^\infty a_k\sin(kt) + b_k\cos(kt).
\end{align}
We now take $\lambda(t)$ as the rate of arrivals at time $t$ in the Erlang-A model. Under this setting, we derive the following expression for the cumulant moment generating function of the fluid approximation and its corresponding partial differential equation whenever the arrival rate is sufficiently large. We do so through a series of technical lemmas. First, we bound the fluid mean when the arrival rate and initial value are sufficiently large.
\begin{lemma}\label{cLemmaNS}
Suppose that $\underline{\lambda} \equiv \inf_{t \geq 0} \lambda(t) > c\mu$ and that $\E{Q_f(0)} > c$. Then,
$$
\E{Q_f(t)} > c
$$
for all time $t \geq 0$.
\begin{proof}
We have seen that $\E{Q_f(t)}$ evolves according to
  $$
  \updot{\mathrm{E}}\left[Q_f(t)\right] = \lambda(t) - \mu (\E{Q_f(t)} \wedge c) - \theta(\E{Q_f(t)} - c)^+
  $$
  at all times $t$. Now, suppose that $\hat t > 0$ is a time such that $\E{Q_f(\hat t)} = c + \epsilon$ for some $\epsilon > 0$. Then, if $\epsilon < \frac{\underline{\lambda} - c\mu}{\theta}$ we have that
  $$
  \updot{\mathrm{E}}\left[Q_f(\hat t)\right] = \lambda(\hat t) - c \mu  - \theta\epsilon \geq \underline{\lambda} - c \mu - \theta \epsilon > 0 .
  $$
  By the continuity of the fluid mean and the fact that $\E{Q_f(0)} = q(0) > c$, we see that $\E{Q_f(t)} > c$ for all time $t \geq 0$.
\end{proof}
\end{lemma}
With this in hand, we now also provide the moment generating function for an $M/M/\infty$ queue with nonstationary arrival rate $\lambda(t)$, which we will use for comparison later in this section.

\begin{lemma}\label{mminfMGF}
Let $Q_\infty(t)$ be an infinite server queue with nonstationary Poisson arrival rate $\lambda(t)$ and exponential service rate $\mu$ and initial value $Q_\infty(t) = q_0$. Then,
$$
\E{e^{\alpha Q_\infty(t)}} =
e^{(e^{\alpha} - 1)\left(\frac{\lambda_0}{\mu}(1 - e^{- \mu t}) +  \sum_{k=1}^\infty \frac{(a_k\mu + b_k k)\sin(k t) + (b_k\mu - a_k k)(\cos(k t) - e^{- \mu t}) }{\mu^2 + k^2}\right)}
\left(e^{- \mu t}(e^{\alpha} - 1) + 1\right)^{q_0}
$$
for all $t \geq 0$ and $\alpha \in \mathbb{R}$.
\begin{proof}
To start, we have that time derivative of the MGF is
$$
\frac{\mathrm{d}\E{e^{\alpha Q_\infty(t)}}}{\mathrm{d} t} = \lambda(t)(e^\alpha - 1)\E{e^{\alpha Q_\infty(t)}} + \mu(e^{-\alpha} - 1)\E{Q_\infty(t) e^{\alpha Q_\infty(t)}}
$$
where $\lambda(t)$ is as defined previously:
$$
\lambda(t) = \lambda_0 + \sum_{k=1}^\infty a_k\sin(kt) + b_k\cos(kt).
$$
This differential equation can be view as a partial differential equation when expressed as
$$
\mu(1 - e^{-\alpha})\frac{\partial M(\alpha, t)}{\partial \alpha}  + \frac{\partial M(\alpha, t)}{\partial t} = \lambda(t)(e^{\alpha} - 1)M(\alpha, t)
$$
where $M(\alpha, t)$ is the moment generating function at time $t$ and space parameter $\alpha$. To simplify our effort, we instead consider the differential equation for the cumulant MGF, which is $G(\alpha , t) = \log(M(\alpha, t))$. This PDE is
$$
\mu(1 - e^{-\alpha})\frac{\partial G(\alpha, t)}{\partial \alpha}  + \frac{\partial G(\alpha, t)}{\partial t} = \lambda(t)(e^{\alpha} - 1)
$$
with the initial condition that
$$
G(\alpha, 0) = \log\left(\E{e^{\alpha Q_\infty(0)}}\right) = \log\left(e^{\alpha q_0}\right) = \alpha q_0.
$$
Using the notation that $G_x = \frac{\partial G}{\partial x}$, we seek to solve the system
$$
\begin{cases}
\mu(1 - e^{-\alpha})G_\alpha + G_t = \lambda(t)(e^{\alpha} - 1)\\
G(\alpha, 0) = \alpha q_0
\end{cases}
$$
and we do so via the method of characteristics. For this approach we introduce variables the characteristic variables $r$ and $s$ and establish the characteristic equations, which are ODE's, as
\begin{align*}
\frac{\mathrm{d}\alpha}{\mathrm{d}s}(r,s) &= \mu(1 - e^{-\alpha}),\\
\frac{\mathrm{d}t}{\mathrm{d}s}(r,s) &= 1,\\
\frac{\mathrm{d}g}{\mathrm{d}s}(r,s) &= \lambda(t)(e^{\alpha} - 1)
\end{align*}
with the initial conditions
\begin{align*}
\alpha(r,0) &= r,\\
t(r,0) &= 0,\\
g(r,0) &= rq_0.
\end{align*}
We can first see that the ODE's for $\alpha$ and $t$ solve to
\begin{align*}
\alpha(r,s) = \log(e^{c_1(r) + \mu s} + 1) &\longrightarrow \alpha(r,s) = \log\left((e^r - 1) e^{\mu s} + 1\right)\\
t(r,s) = s + c_2(r) &\longrightarrow t(r,s) = s
\end{align*}
and so we can now use these to solve the remaining ODE. After substituting we have
$$
\frac{\mathrm{d}g}{\mathrm{d}s}(r,s) = \lambda(s)(e^r - 1) e^{\mu s}
$$
which gives a solution of
\begin{align*}
g(r,s) &= (e^r - 1)\left(\frac{\lambda_0}{\mu}(e^{\mu s} - 1) +  \sum_{k=1}^\infty \frac{(a_k\mu + b_k k)\sin(k s)e^{\mu s} + (b_k\mu - a_k k)(\cos(k s)e^{\mu s} - 1) }{\mu^2 + k^2}\right) + rq_0.
\end{align*}
So, using $s = t$ and $r = \log\left(e^{- \mu t}(e^{\alpha} - 1) + 1\right)$, we have that
\begin{align*}
G(\alpha, t)
&=
g(\log\left(e^{- \mu t}(e^{\alpha} - 1) + 1\right), t)
\\&=
(e^{\alpha} - 1)\left(\frac{\lambda_0}{\mu}(1 - e^{- \mu t}) +  \sum_{k=1}^\infty \frac{(a_k\mu + b_k k)\sin(k t) + (b_k\mu - a_k k)(\cos(k t) - e^{- \mu t}) }{\mu^2 + k^2}\right)
\\&\qquad+
\log\left(e^{- \mu t}(e^{\alpha} - 1) + 1\right)q_0
\end{align*}
and therefore by solving for $M(\alpha, t) = e^{G(\alpha, t)}$ we attain the stated result.
\end{proof}
\end{lemma}

Now that we have established these lemmas we proceed with the analysis of the nonstationary Erlang-A. In the next theorem we give explicit forms for the fluid form of the cumulant MGF and its corresponding partial differential equation.

\begin{theorem}\label{fluidNSMGF}
If $ \inf_{t \leq \infty} \lambda(t) \equiv \underline{\lambda} > c\mu$ and $q(0) > c$, then in for all $t \geq 0$ we have that
     \begin{align}
    \frac{\partial G_f(t,\alpha) }{\partial t} &=  \lambda(t) \cdot (e^{\alpha} -1) +  \theta \cdot (e^{-\alpha} -1) \cdot \frac{\partial G_f(t,\alpha)}{\partial \alpha} - c \cdot ( \theta - \mu) \cdot (e^{-\alpha} -1) \label{fluidNSMGFeq1}
  \end{align}
which gives a solution of
     \begin{align}
   G_f(t, \alpha)
  &=
  (e^\alpha - 1)\left(\frac{\lambda_0}{\theta}(1 - e^{-\theta t}) +  \sum_{k=1}^\infty \frac{(a_k\theta + b_k k)\sin(k t) + (b_k\theta - a_k k)(\cos(k t) - e^{-\theta t}) }{\theta^2 + k^2}\right)
  \nonumber
  \\&\qquad
  +
  \frac{c(\theta - \mu)}{\theta}\alpha
  +
  \log((e^\alpha - 1)e^{-\theta t} + 1)\left( q(0) - \frac{c(\theta - \mu)}{\theta} \right)
  \end{align}
  for all $t \geq 0$ and all $\alpha \in \mathbb{R}$.
  \begin{proof}
 From Equation~\ref{cumulantMGF}, we have that the PDE for the fluid approximation's cumulant moment generating function is
  $$
  \frac{\partial G_f(t,\alpha) }{\partial t} =  \lambda(t) (e^{\alpha} -1) +  \theta  (e^{-\alpha} -1) \frac{\partial G_f(t,\alpha)}{\partial \alpha}
     - ( \theta - \mu) (e^{-\alpha} -1)    \left( \frac{\partial G(t,\alpha) }{\partial \alpha} \wedge c  \right) .
  $$
  Now, recall that $\frac{\partial G_f(t,\alpha) }{\partial \alpha} = \frac{\E{Q_f(t)e^{\alpha Q_f(t)}}}{\E{e^{\alpha Q_f(t)}}}$. Using the FKG inequality and our observation from Lemma \ref{cLemmaNS} that $\E{Q_f(t)} > c$, we have that
  $$
  \E{Q_f(t)e^{\alpha Q_f(t)}} \geq \E{Q_f(t)}\E{e^{\alpha Q_f(t)}} > c \E{e^{\alpha Q_f(t)}}
  $$
  and so $\left(\frac{\partial G_f(t,\alpha) }{\partial \alpha} \wedge c\right) = c$. Thus, we have the PDE given in Equation~\ref{fluidNSMGFeq1} and so now we seek to find it's solution. We approach this via the method of characteristics. Because $G_f(0, \alpha) = \log(\E{e^{\alpha Q_f(0)}}) = \alpha q(0)$, we see that we seek to solve the following system
  $$
  \begin{cases}
  \theta(1 - e^{-\alpha})G_{(\alpha)} + G_{(t)} = \lambda(t) (e^{\alpha} - 1) + c(\theta - \mu)(1 - e^{-\alpha})\\
  G_f(0,\alpha) = \alpha q(0)
  \end{cases}
  $$
  where $G_{(x)} = \frac{\partial G_f}{\partial x}$. Introducing characteristic variables $r$ and $s$, we have the characteristic ODE's as
  \begin{align*}
  \frac{\mathrm{d}\alpha}{\mathrm{d}s}(r,s) &= \theta (1 - e^{-\alpha})\\
  \frac{\mathrm{d}t}{\mathrm{d}s}(r,s) &= 1\\
  \frac{\mathrm{d}g}{\mathrm{d}s}(r,s) &= \lambda(t)(e^{\alpha} - 1) + c(\theta - \mu)(1-e^{-\alpha})
  \end{align*}
  with initial conditions $\alpha(r, 0) = r$, $t(r,0) = t$, and $g(r,0) = rq(0)$. Then, we can solve the first two ODE's to see that
  \begin{align*}
  \alpha(r,s) &= \log((e^r - 1)e^{\theta s} + 1)\\
  t(r,s) &= s
  \end{align*}
  and so we can use these to solve the remaining equation. Substituting in, we have the ODE as
  $$
  \frac{\mathrm{d}g}{\mathrm{d}s}(r,s) = \lambda(s)e^{\theta s}(e^r - 1) + c(\theta - \mu)\frac{e^{\theta s}(e^r - 1)}{e^{\theta s}(e^r - 1) + 1}
  $$
  and this now solves to
  \begin{align*}
  g(r,s)
  &=
  (e^r - 1)\left(\frac{\lambda_0}{\theta}(e^{\theta s} - 1) +  \sum_{k=1}^\infty \frac{(a_k\theta + b_k k)\sin(k s)e^{\theta s} + (b_k\theta - a_k k)(\cos(k s)e^{\theta s} - 1) }{\theta^2 + k^2}\right)
  \\&\qquad
  +
  \frac{c(\theta - \mu)}{\theta}\left(\log\left((e^r-1)e^{\theta s} + 1\right) - r\right)
  +
  r q(0) .
  \end{align*}
  Now, we can rearrange our solutions to find $s = t$ and $r = \log((e^\alpha - 1)e^{-\theta t} + 1)$. Then, we have that
  \begin{align*}
  G_f(t, \alpha)
  &=
  g(\log((e^\alpha - 1)e^{-\theta t} + 1), t)
  \\
  &=
  (e^\alpha - 1)e^{-\theta t}\left(\frac{\lambda_0}{\theta}(e^{\theta t} - 1) +  \sum_{k=1}^\infty \frac{(a_k\theta + b_k k)\sin(k t)e^{\theta t} + (b_k\theta - a_k k)(\cos(k t)e^{\theta t} - 1) }{\theta^2 + k^2}\right)
  \\&\qquad
  +
  \frac{c(\theta - \mu)}{\theta}\left(\alpha - \log((e^\alpha - 1)e^{-\theta t} + 1)\right)
  +
  \log((e^\alpha - 1)e^{-\theta t} + 1) q(0)
  \end{align*}
  and this simplifies to the stated result.
  \end{proof}
\end{theorem}

Like the approach in our investigation of the steady-state scenario, we can now observe that the fluid approximation is equivalent in distribution to a infinite server queue shifted by $\gamma \equiv \frac{c(\theta - \mu)}{\theta}$. This gives rise to the following.

\begin{theorem}
For the Erlang-A queue with nonstationary arrival rate $\lambda(t)$ such that $\underline{\lambda} \equiv \inf_{t \geq 0} \lambda(t) > c\mu$ and initial value $q(0) > c$,  the fluid approximation of the MGF is equivalent to that of a shifted $M/M/\infty$ queue with arrival rate $\lambda(t)$, service rate $\theta$, initial value $q(0) - \frac{c(\theta - \mu)}{\theta}$, and linear shift $\frac{c(\theta - \mu)}{\theta}$.
\begin{proof}
Observe from Theorem~\ref{fluidNSMGF} that the fluid MGF for the Erlang-A under these conditions is
\begin{align*}
&M_f(t, \alpha)
=
e^{G_f(t , \alpha)}
\\
&=
e^{
(e^\alpha - 1)\left(\frac{\lambda_0}{\theta}(1 - e^{-\theta t}) +  \sum_{k=1}^\infty \frac{(a_k\theta + b_k k)\sin(k t) + (b_k\theta - a_k k)(\cos(k t) - e^{-\theta t}) }{\theta^2 + k^2}\right)
+
  \frac{c(\theta - \mu)}{\theta}\alpha
}
\left((e^\alpha - 1)e^{-\theta t} + 1\right)^{ q(0) - \frac{c(\theta - \mu)}{\theta} }
\end{align*}
which is of a form that we can recognize. Comparing it to Lemma~\ref{mminfMGF}, we can see that $Q_f$ is of the form of a shifted $M/M/\infty$ queue with arrival rate $\lambda(t)$, service rate $\theta$, initial value $q(0) - \frac{c(\theta - \mu)}{\theta}$, and linear shift $\frac{c(\theta - \mu)}{\theta}$, thus enforcing that the fluid model does start at $q(0)$.
\end{proof}
\end{theorem}

This representation of the fluid approximation allows us to now provide upper and lower bounds for the moments of the Erlang-A system.

\begin{corollary}
Let $Q(t)$ represent the Erlang-A queue with nonstationary arrival rate $\lambda(t)$ such that $\underline{\lambda} \equiv \inf_{t \geq 0} \lambda(t) > c\mu$ and initial value $q(0) > c$, and let $Q_f(t)$ represent the corresponding fluid approximation. Then, if $\theta > \mu$
$$
\E{\left(Q_f(t) - \gamma\right)^m} \leq \E{Q(t)^m} \leq \E{Q_f(t)^m}
$$
 and if $\theta < \mu$
$$
\E{Q_f(t)^m} \leq \E{Q(t)^m} \leq \E{\left(Q_f(t) - \gamma\right)^m}
$$
for all time $t > 0$ and all $m \in \mathbb{Z}^+$, where $\gamma =  \frac{c(\theta - \mu)}{\theta}$.
\begin{proof}
In each case, the bound involving the fluid approximation of the moment is a direct consequence of Theorem~\ref{mMoment} and so only the other two bounds remain to be shown. We now note that since we have characterized the fluid approximation as a shifted $M/M/\infty$ queue, the remaining bounds are from the unshifted version of this system and, by following the same arguments as in Theorems~\ref{sandwichthm} and~\ref{sandwichthm2} regarding the rates of departure in the corresponding states of the Erlang-A queue and the $M/M/\infty$ queue, this completes the proof.
\end{proof}
\end{corollary}

\subsection{Numerical Results}\label{mgfnumres}

\begin{figure}[H]
\centering
\hspace{-.35in}~\includegraphics[width = .5\textwidth]{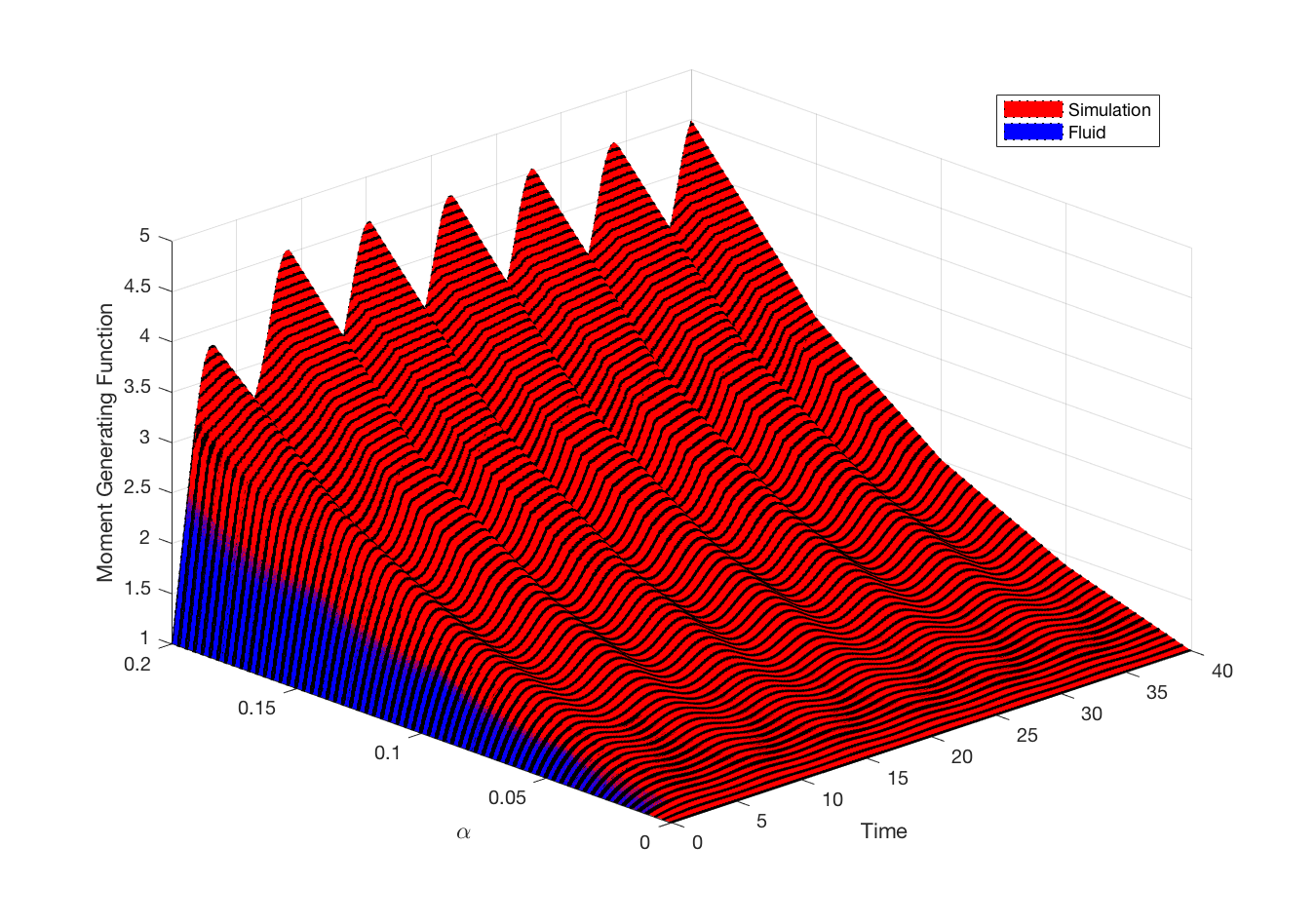}~\hspace{-.15in}~\includegraphics[width = .5\textwidth]{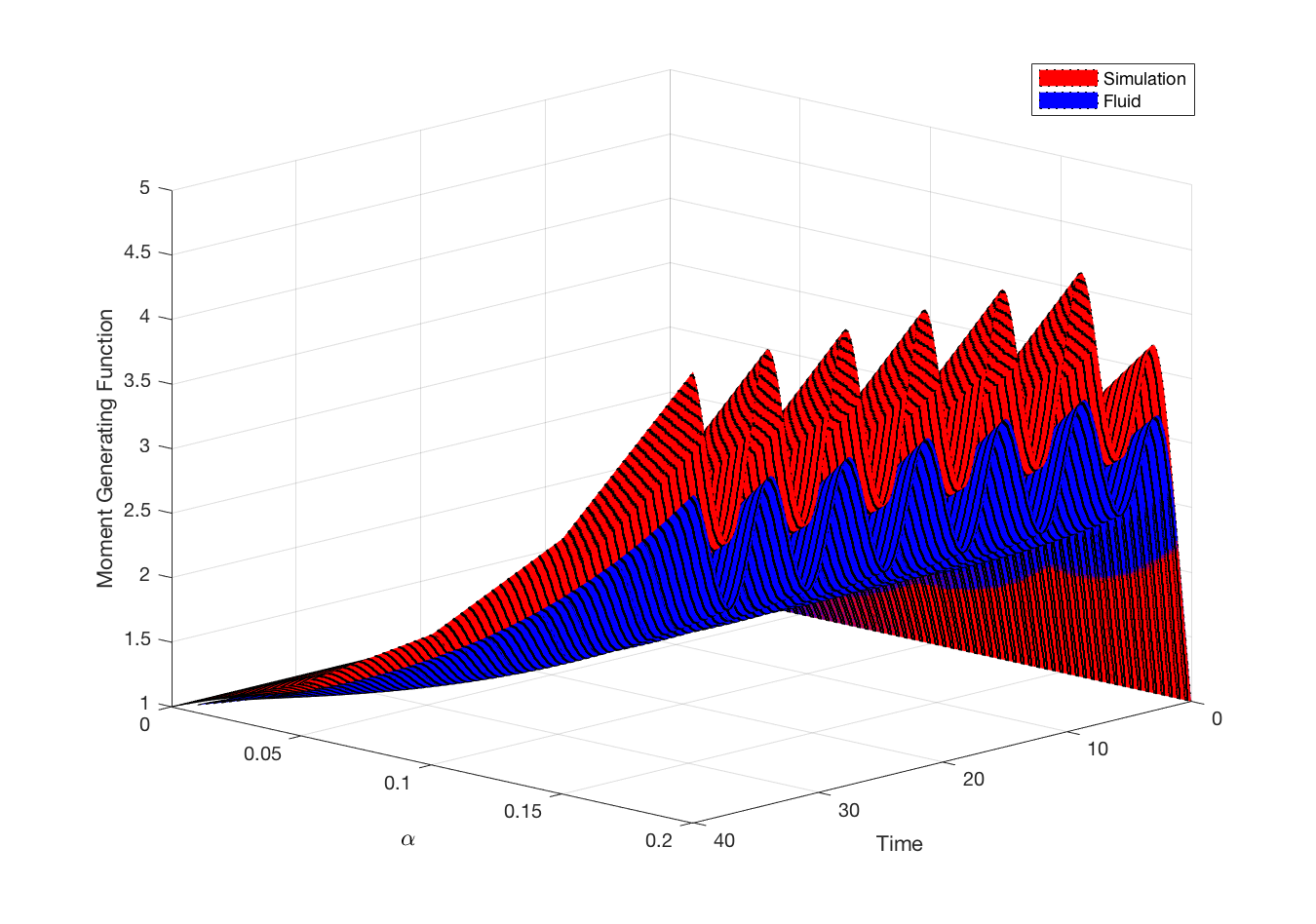}
\captionsetup{justification=centering}
 \caption{$\lambda(t) = 5 + \sin (t) $, $\mu = 1$, $\theta = 0.5$, $Q(0) = 0$, $c=5$. \\ Front view (left) and rear view (right). \label{MGFFig3} }
\end{figure}

\begin{figure}[H]
\centering
\hspace{-.35in}~\includegraphics[width = .5\textwidth]{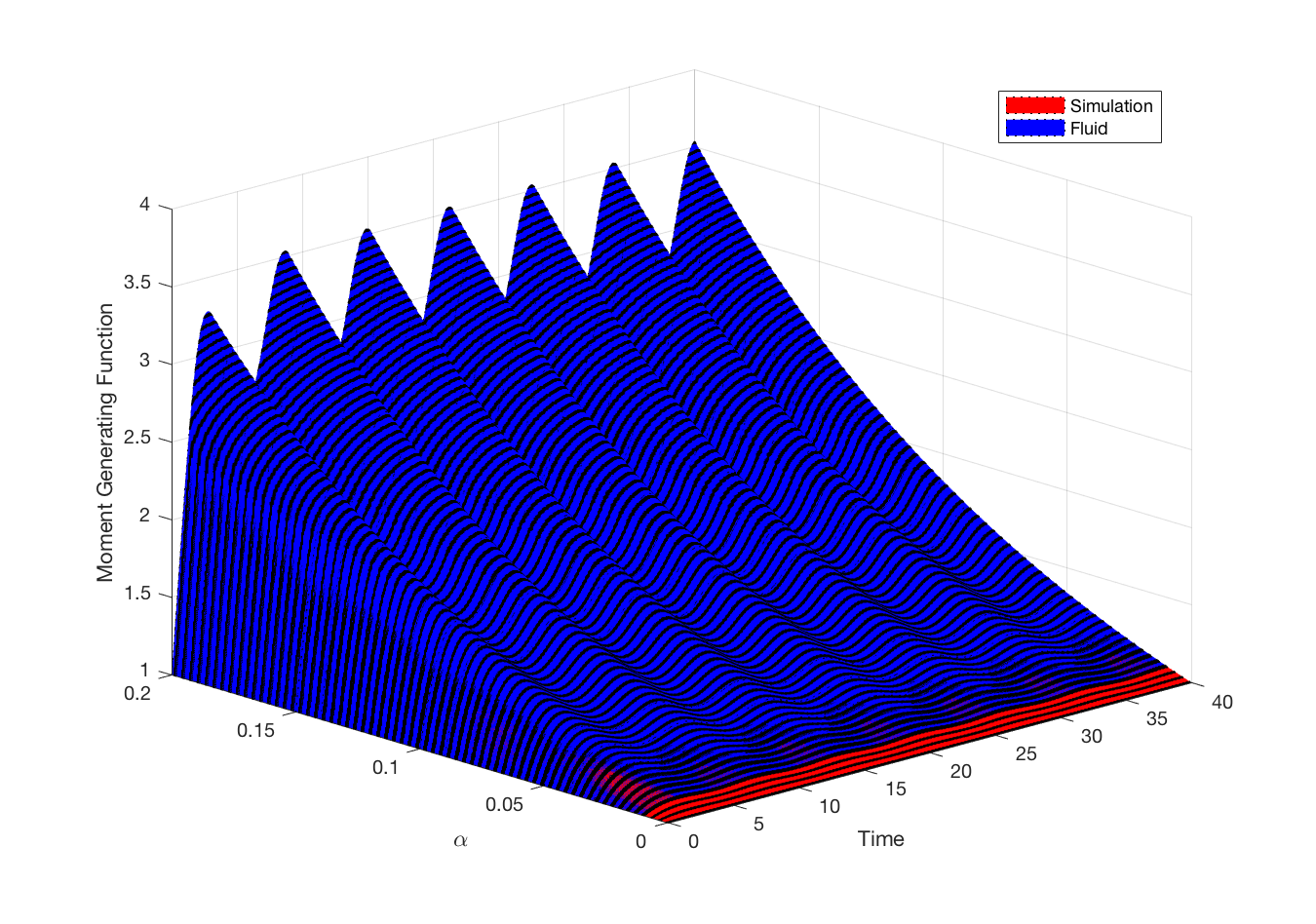}~\hspace{-.15in}~\includegraphics[width = .5\textwidth]{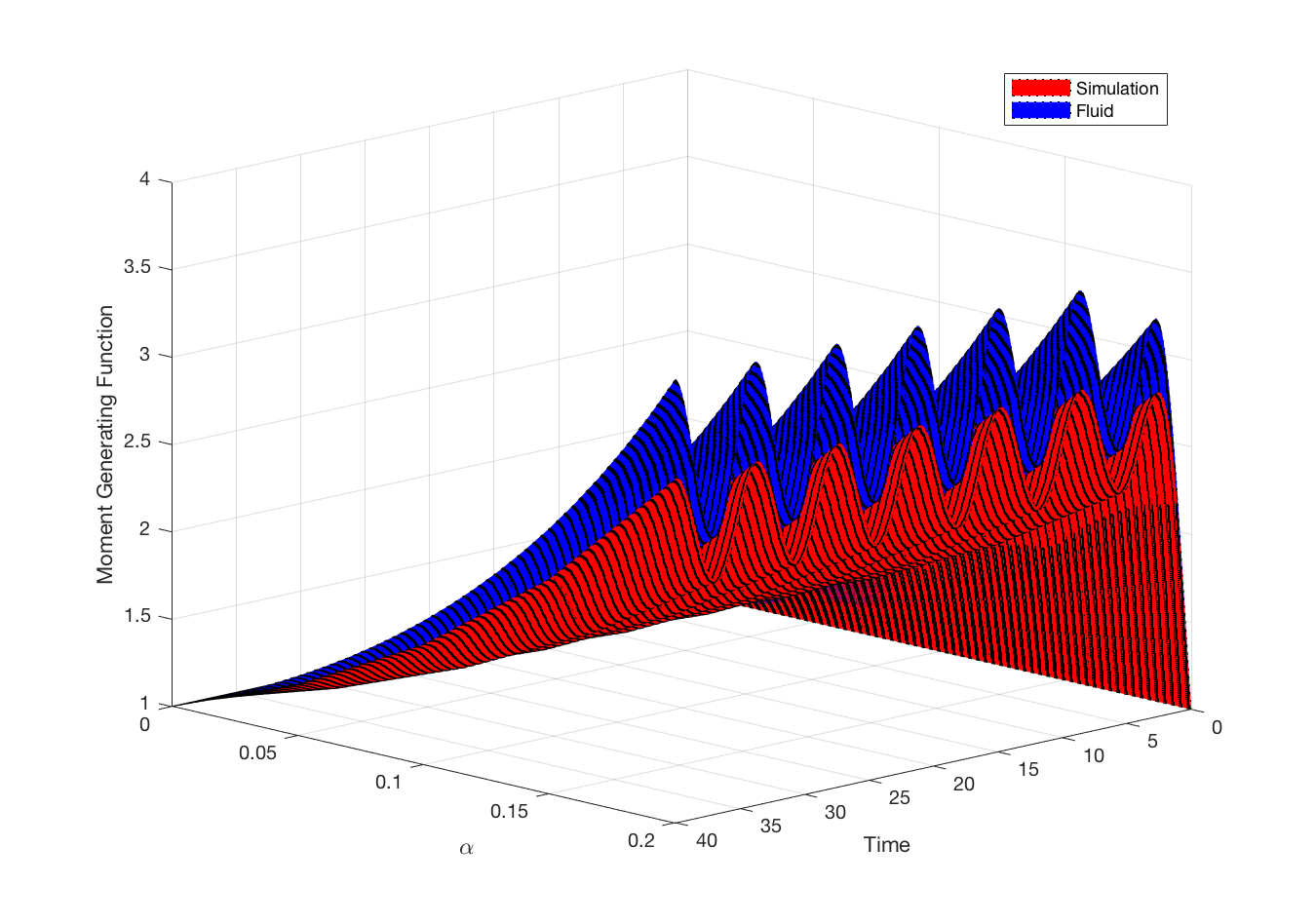}
\captionsetup{justification=centering}
 \caption{$\lambda(t) = 5 +  \sin (t) $, $\mu = 1$, $\theta = 2$, $Q(0) = 0$, $c=5$. \\ Front view (left) and rear view (right). \label{MGFFig4} }
\end{figure}

\begin{figure}[H]
\centering
\hspace{-.35in}~\includegraphics[width = .5\textwidth]{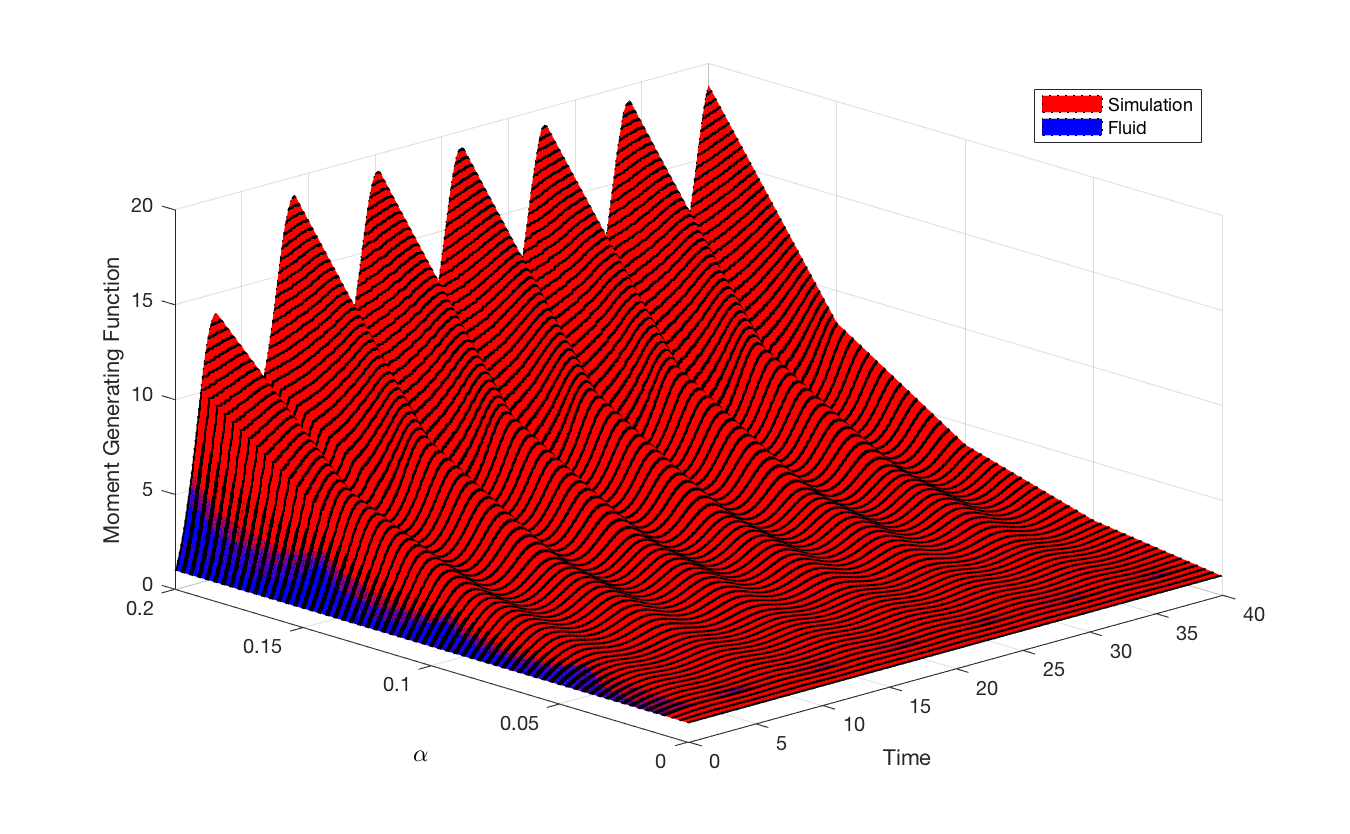}~\hspace{-.15in}~\includegraphics[width = .5\textwidth]{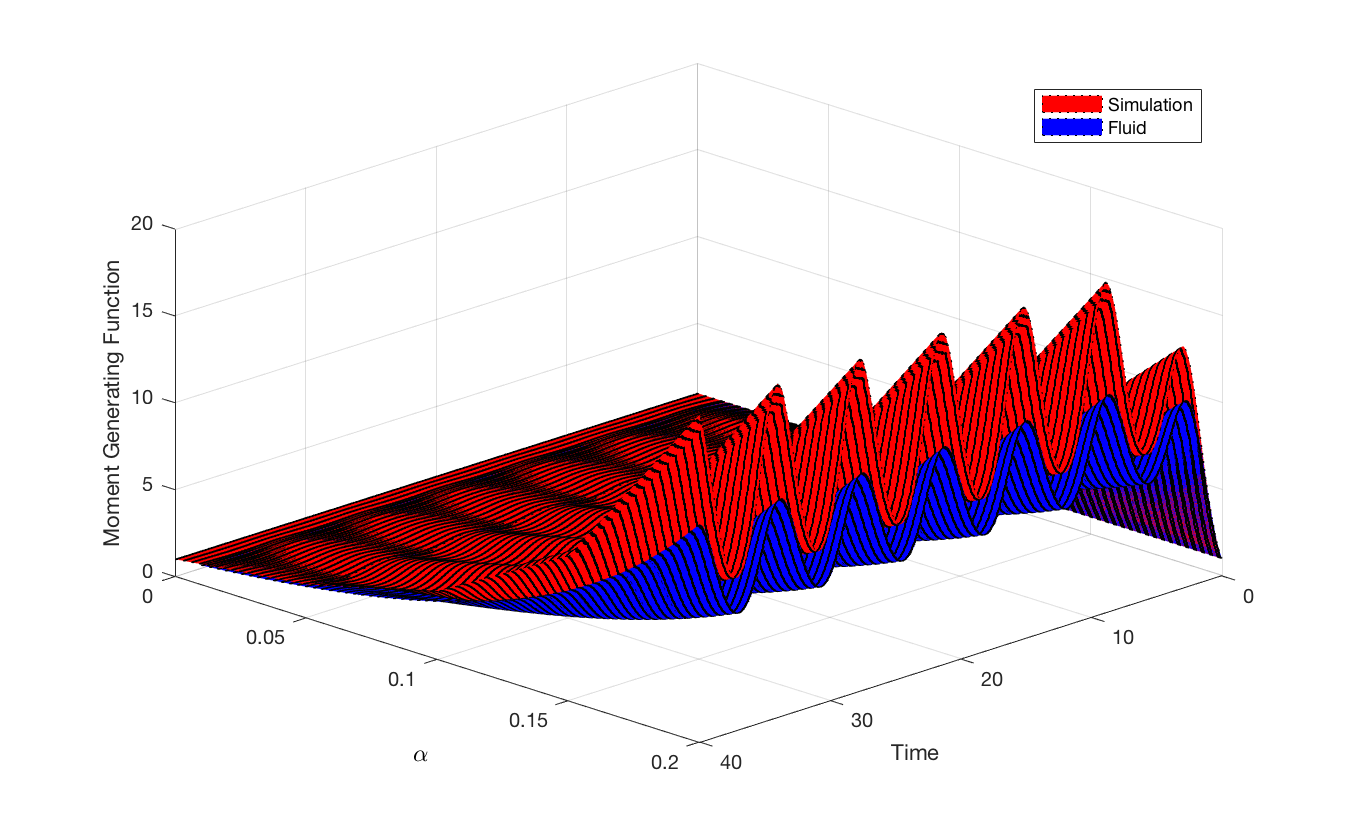}
\captionsetup{justification=centering}
 \caption{$\lambda(t) = 10 + 2 \cdot \sin (t) $, $\mu = 1$, $\theta = 0.5$, $Q(0) = 0$, $c=10$. \\ Front view (left) and rear view (right). \label{MGFFig1} }
\end{figure}

\begin{figure}[H]
\centering
\hspace{-.35in}~\includegraphics[width = .5\textwidth]{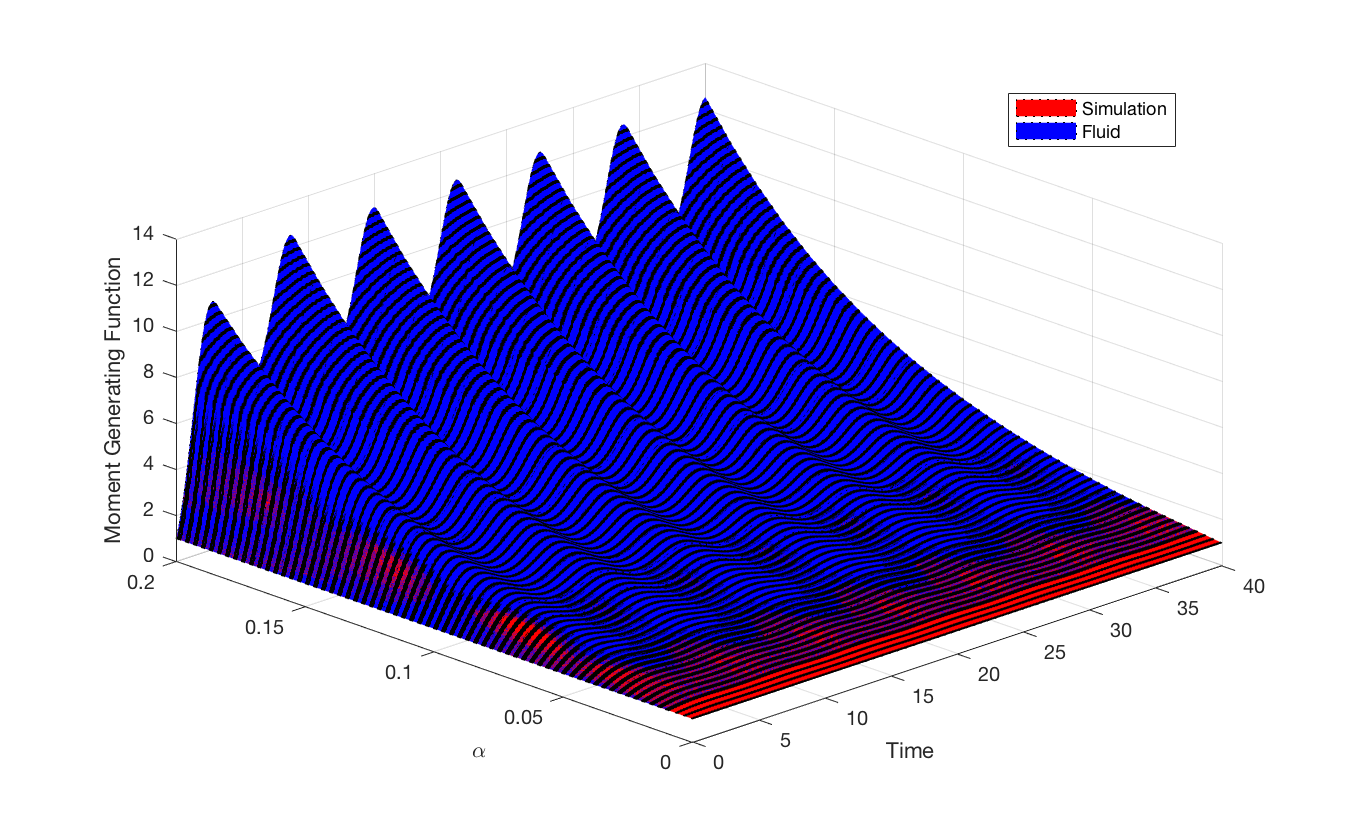}~\hspace{-.15in}~\includegraphics[width = .5\textwidth]{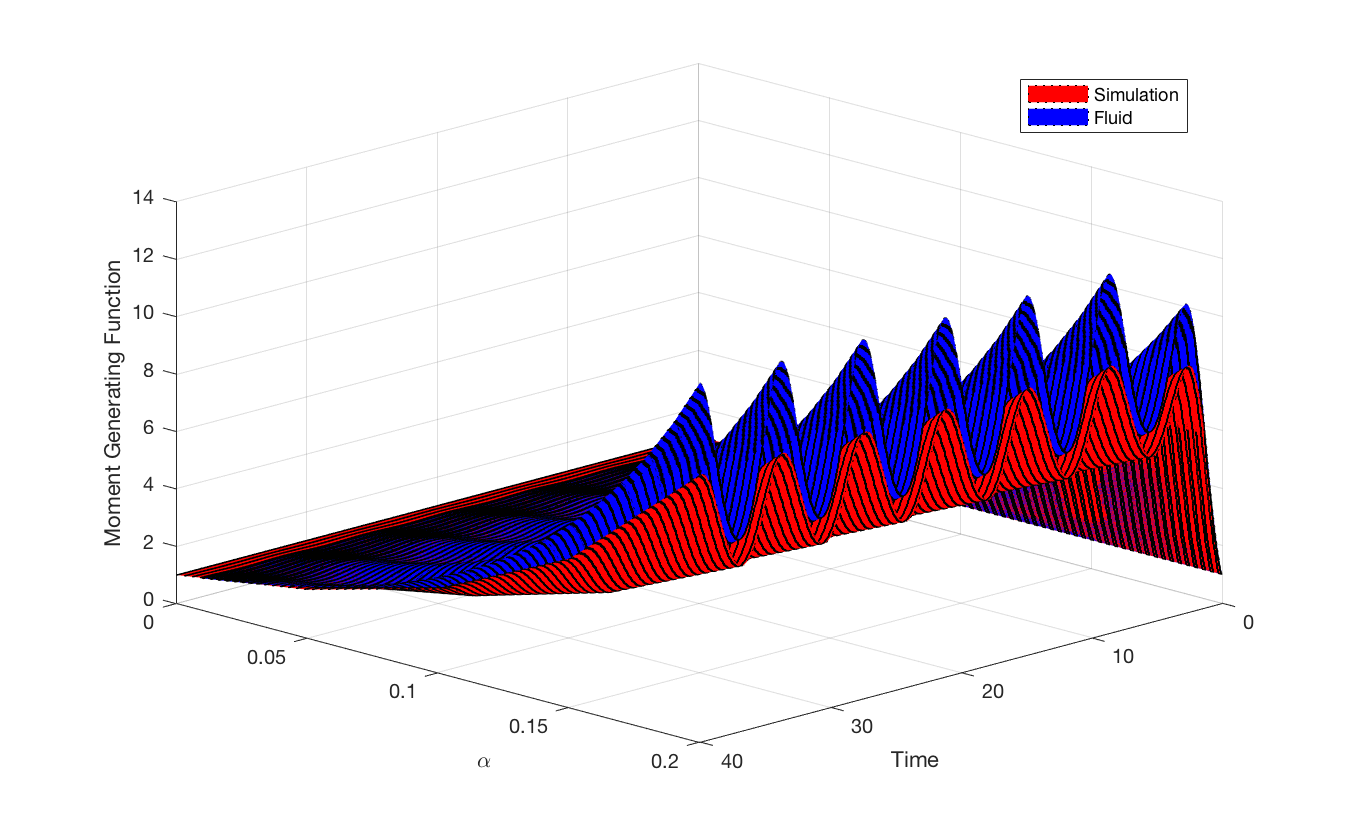}
\captionsetup{justification=centering}
 \caption{$\lambda(t) = 10 + 2 \cdot \sin (t) $, $\mu = 1$, $\theta = 2$, $Q(0) = 0$, $c=10$. \\ Front view (left) and rear view (right). \label{MGFFig2} }
\end{figure}

In this subsection we describe various numerical experiments demonstrating these findings. We first have Figures~\ref{MGFFig3},~\ref{MGFFig4},~\ref{MGFFig1}, and~\ref{MGFFig2}, which compare simulated value of the moment generating function to their fluid approximations. In the first two figures, the arrival intensity is $\lambda(t) = 5 + \sin(t)$, the service rate is $\mu = 1$, and the number of servers is $c = 5$. The abandonment rates are the differing component of these plots, with $\theta = 0.5$ and $\theta = 2$ as the two respective values.  These same comparisons are made in the latter two figures, however in this case the arrival rate is instead $\lambda(t) = 10 + 2\sin(t)$ and the number of servers is $c = 10$. \\

Through these plots one can observe that the true MGF dominates the fluid approximation when $\theta < \mu$ and that the fluid dominates the stochastic value when $\theta > \mu$. This is of course stated with the understanding that for small values of $\alpha$ or for times near 0 the values of the MGF and the approximation are quite close and so with numerical error the surfaces may overlap.

\begin{figure}[H]
        \begin{subfigure}[b]{0.475\textwidth}
            \centering
            \includegraphics[scale = .45]{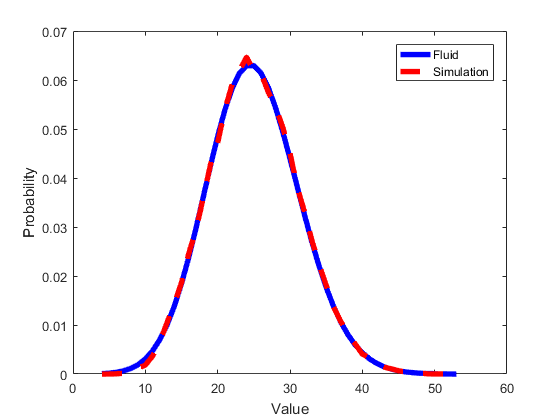}
            \caption{{\small $\theta = 0.5$, $c = 15$}}
            \label{fig:3}
        \end{subfigure}
        \quad
        \begin{subfigure}[b]{0.475\textwidth}
            \centering
            \includegraphics[scale = .45]{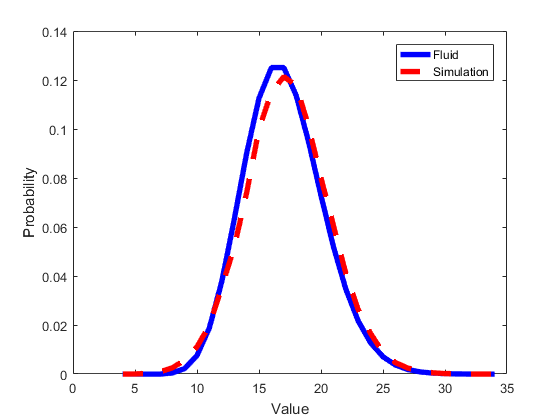}
            \caption{{\small $\theta = 2$, $c = 15$}}
            \label{fig:4}
        \end{subfigure}
        \begin{subfigure}[b]{0.475\textwidth}
            \centering
            \includegraphics[scale = .45]{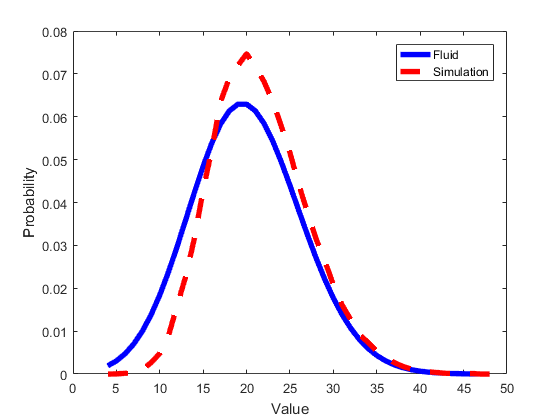}
            \caption{{\small $\theta = 0.5$, $c = 20$}}
            \label{fig:1}
        \end{subfigure}
        \hfill
        \begin{subfigure}[b]{0.475\textwidth}
            \centering
            \includegraphics[scale = .45]{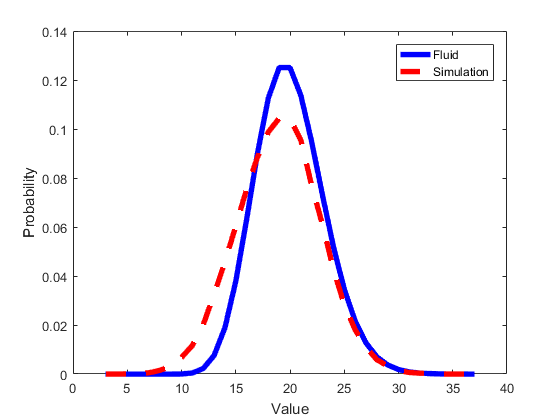}
            \caption{{\small $\theta = 2$, $c = 20$}}
            \label{fig:2}
        \end{subfigure}
        \begin{subfigure}[b]{0.475\textwidth}
            \centering
            \includegraphics[scale = .45]{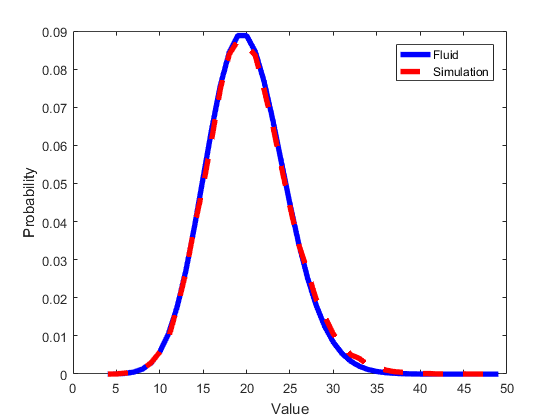}
            \caption{{\small $\theta = 0.5$, $c = 25$}}
            \label{fig:3}
        \end{subfigure}
        \quad
        \begin{subfigure}[b]{0.475\textwidth}
            \centering
            \includegraphics[scale = .45]{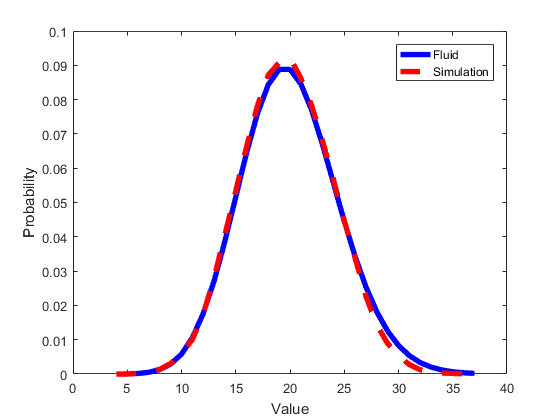}
            \caption{{\small $\theta = 2$, $c = 25$}}
            \label{fig:4}
        \end{subfigure}
        \caption
        {Empirical and Fluid Limiting Distributions for $\lambda = 20$ and $\mu = 1$.}
        \label{limdists}
    \end{figure}

In Figure~\ref{limdists} we plot the limiting distribution for the steady-state Erlang-A. For these plots we take $\lambda = 20$ and $\mu = 1$, and then vary $\theta$ and $c$. For the three plots on the left we take the abandonment rate to be $\theta = 0.5$ and for those on the right we set $\theta = 2$. For the top two plots we set the number of servers as $c = 15$, in the middle two $c = 20$, and in the bottom two we make $c = 25$. We observe that the approximate distribution is quite close when $\lambda$ is not near $c\mu$ but the approximation is less accurate when $\lambda = c\mu$.  This finding is consistent with much of the literature that focuses on finding novel approximations for queueing networks and optimal control of these networks, see for example \citet{hampshire2010dynamic, hampshire2009time, hampshire2009dynamic, pender2017approximations, niyirora2016optimal, qin2017dynamic}. We note here that these approximations are not all of the same form: recall that when $\lambda \geq c\mu$ the fluid approximation is equivalent in distribution to a shifted Poisson random variable with parameter $\frac{\lambda}{\theta}$, but when $\lambda < c\mu$ it is equivalent to a Poisson distribution with parameter $\frac{\lambda}{\mu}$.
\vspace{-.1in}

\begin{figure}[H]
        \begin{subfigure}[b]{0.475\textwidth}
            \centering
            \includegraphics[scale = .45]{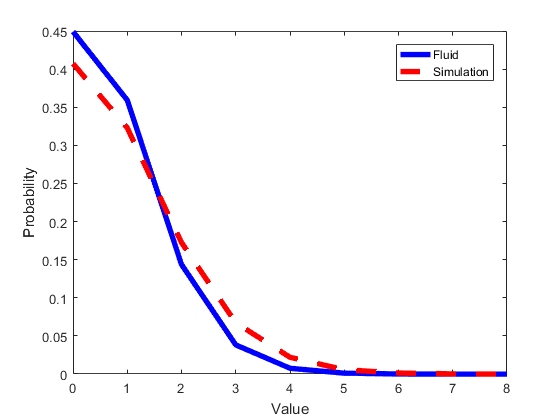}
            \caption{{\small $\theta = 0.5$, $\lambda = 0.8$}}
            \label{fig:3}
        \end{subfigure}
        \quad
        \begin{subfigure}[b]{0.475\textwidth}
            \centering
            \includegraphics[scale = .45]{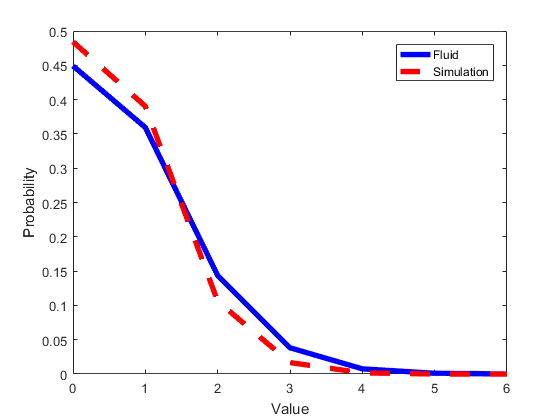}
            \caption{{\small $\theta = 2$, $\lambda = 0.8$}}
            \label{fig:4}
        \end{subfigure}
        \begin{subfigure}[b]{0.475\textwidth}
            \centering
            \includegraphics[scale = .45]{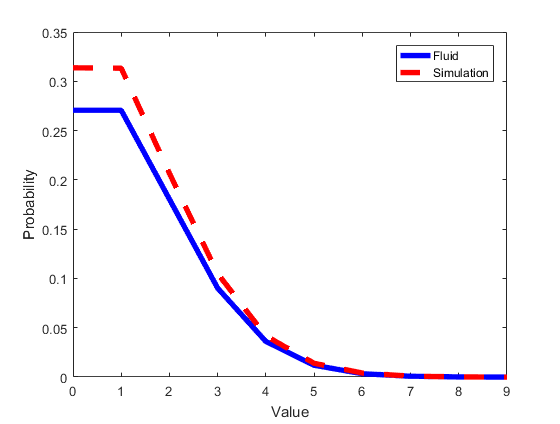}
            \caption{{\small $\theta = 0.5$, $\lambda = 1$}}
            \label{fig:1}
        \end{subfigure}
        \hfill
        \begin{subfigure}[b]{0.475\textwidth}
            \centering
            \includegraphics[scale = .45]{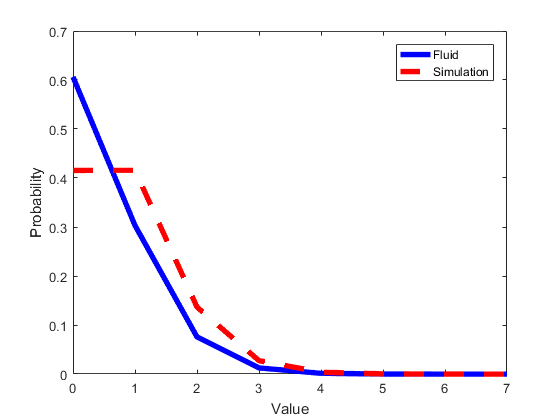}
            \caption{{\small $\theta = 2$, $\lambda = 1$}}
            \label{fig:2}
        \end{subfigure}
        \begin{subfigure}[b]{0.475\textwidth}
            \centering
            \includegraphics[scale = .45]{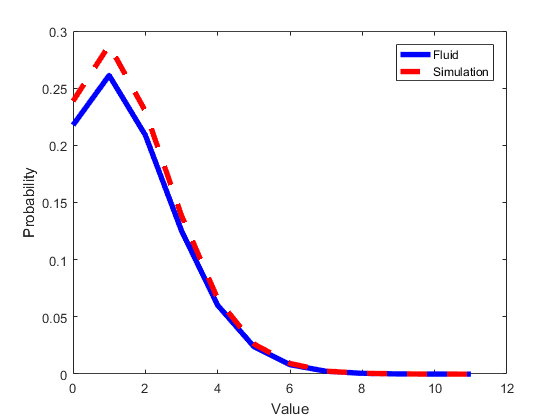}
            \caption{{\small $\theta = 0.5$, $\lambda = 1.2$}}
            \label{fig:3}
        \end{subfigure}
        \quad
        \begin{subfigure}[b]{0.475\textwidth}
            \centering
            \includegraphics[scale = .45]{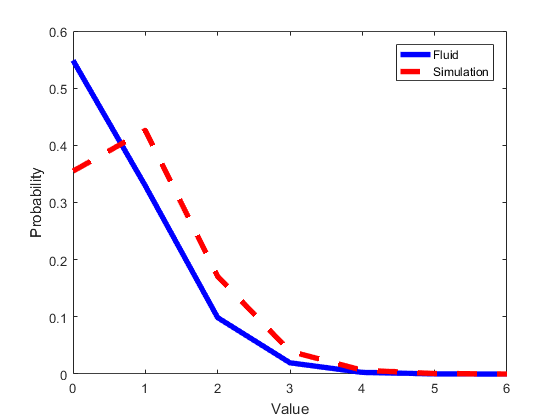}
            \caption{{\small $\theta = 2$, $\lambda = 1.2$}}
            \label{fig:4}
        \end{subfigure}
        \caption
        {Empirical and Fluid Limiting Distributions for $c = 1$ and $\mu = 1$.}
        \label{singlelimdists}
    \end{figure}

    In Figure~\ref{singlelimdists} we examine the limiting distributions for the single server case. In these plots we set $\mu = 1$ and then vary the arrival rate and the abandonment rate. On all plots on the left we set $\theta = 0.5$ and on the right $\theta = 2$. Further, in the top pair we make $\lambda = 0.8$, in the middle we let $\lambda = 1$, and in the bottom pair $\lambda = 1.2$. As in Figure~\ref{limdists}, Figure~\ref{singlelimdists} shows that our approximations are quite good.  Thus, we are able to capture single server dynamics as well as large-scale multi-server dynamics even though they are quite different.  This is even more useful as our approximations are non-asymptotic and don't rely on scaling the number of servers.

\begin{figure}[H]
\centering
\includegraphics[width = .5\textwidth]{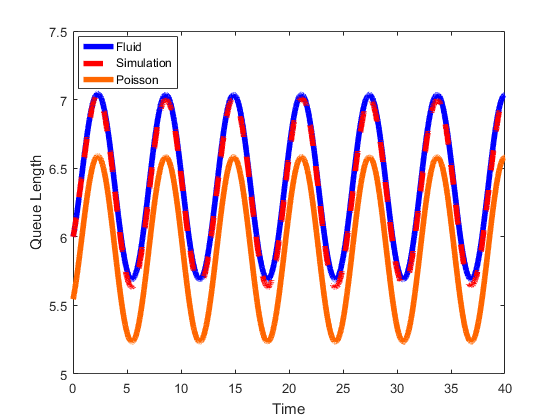}~\hspace{-.15in}~\includegraphics[width = .5\textwidth]{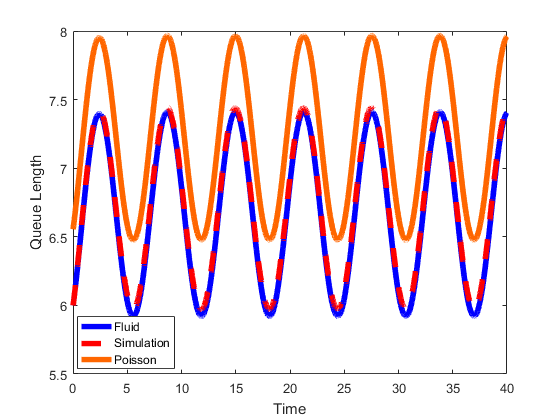}
\captionsetup{justification=centering}
 \caption{Queue Mean for $\lambda(t) = 6.5 + \sin (t) $, $\mu = 1$, $Q(0) = 6$, $c=5$. \\ $\theta = 1.1$ (left) and $\theta = 0.9$ (right). \label{meanpoisfig} }
\end{figure}

\begin{figure}[H]
\centering
\includegraphics[width = .5\textwidth]{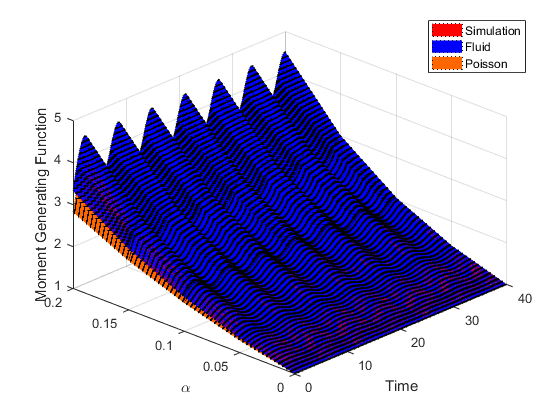}~\hspace{-.15in}~\includegraphics[width = .5\textwidth]{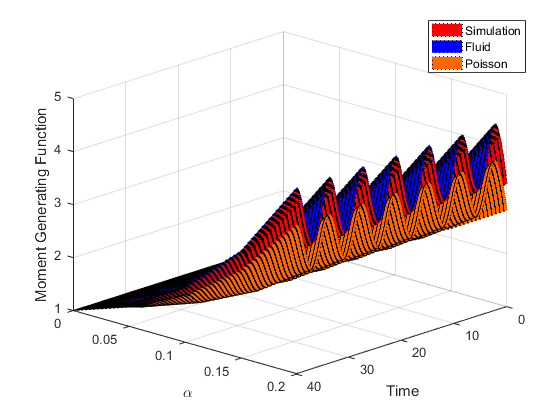}
\captionsetup{justification=centering}
 \caption{MGF for $\lambda(t) = 6.5 + \sin (t) $, $\mu = 1$, $\theta = 1.1$, $Q(0) = 6$, $c=5$. \\ Front view (left) and rear view (right). \label{MGFpoisfig1} }
\end{figure}

\begin{figure}[H]
\centering
\includegraphics[width = .5\textwidth]{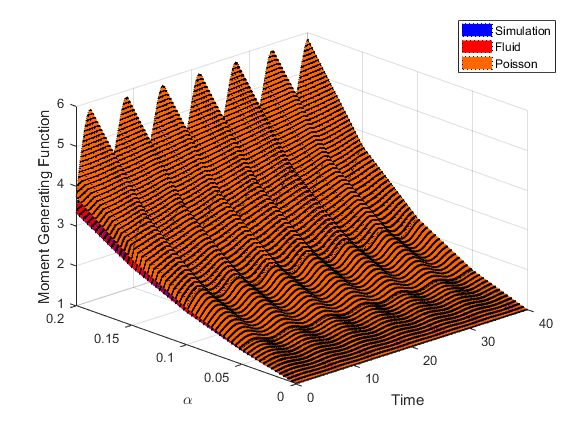}~\hspace{-.15in}~\includegraphics[width = .5\textwidth]{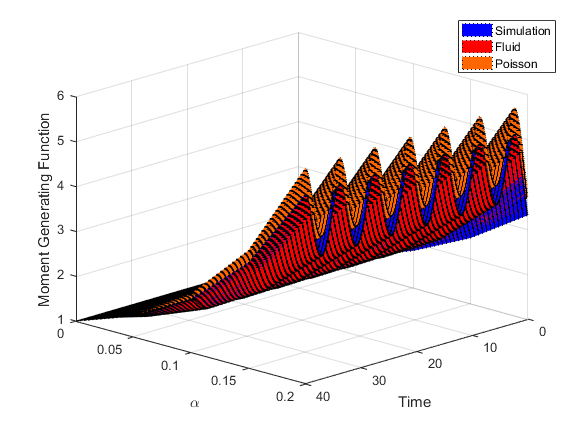}
\captionsetup{justification=centering}
 \caption{MGF for $\lambda(t) = 6.5 + \sin (t) $, $\mu = 1$, $\theta = 0.9$, $Q(0) = 6$, $c=5$. \\ Front view (left) and rear view (right). \label{MGFpoisfig2} }
\end{figure}

In Figures~\ref{meanpoisfig},~\ref{MGFpoisfig1}, and~\ref{MGFpoisfig2}, we take the arrival rate as $\lambda(t) = 6.5 + \sin (t)$, the service rate as $\mu = 1$, and the number of servers as $c = 5$. Because $\inf_{t \geq 0} \lambda(t) > c \mu$, we use the characterization of the fluid approximation as a shifted $M/M/\infty$ queue and compare the simulated system, the fluid approximation, and the unshifted $M/M/\infty$. In the first figure we consider the mean for $\theta = 1.1$ and and $\theta = 0.9$ and find that while the fluid approximation is quite close the unshifted system is not near to the Erlang-A system, even for these relatively similar rates of service and abandonment. We find the same for the latter two figures, in which we plot the moment generating function for $\theta = 1.1$ and $\theta = 0.9$, respectively.

\section{Conclusion}

In this paper we have investigated the Erlang-A queueing system through comparison to the fluid approximations of its moments and moment generating function as well as of its cumulants and cumulant moment generating function. Through recognizing the convexity in the differential equations describing these approximations, we have found fundamental relationships between the values of these quantities and their fluid counterparts: when the rate of abandonment is less than the rate of service the true value dominates the approximation, when the service rate is larger the approximation dominates the true value, and when the rates of abandonment and service are equal, the two are equivalent.

In forming these inequalities, we have found explicit representations of the fluid approximations through equivalences in distribution with Poisson random variables and infinite server queues, in the stationary and non-stationary cases, respectively. These characterizations both give insight into the approximations themselves and yield natural inequalities that complement those from the approximations. We have demonstrated the performance of these bounds through simulations. Through consideration of both these findings and the empirical experiments, we can identify interesting directions of future work.

For example, it would be of great interest to gain more explicit insights into the gap between the fluid approximations and the true values. This is a non-trivial endeavor, which stems from the non-differentiablility and non-closure in the differential equations for the true expectations. The numerical experiments in this work indicate that the fluid approximations may often be quite close but not exact, and additional understanding would be useful in practice.  Moreover, extending our results to more complicated queueing systems where the arrival and service processes follow phase type distributions is of interest given the new work of \citet{pender2017approximations, ko2016strong, ko2017diffusion}.

Additionally, it would be even more useful to gain a better understanding of the limiting distribution of the Erlang-A queue. As we discuss in the paper, the empirical experiments in Subsection~\ref{mgfnumres} indicate that the true limiting distributions closely resemble the shifted Poisson distributions that we have found as characterizations of our fluid approximations. In particular, the approximations seem quite close when $\lambda$ is not near $c\mu$. As a simple extension of this work, it can be observed that some sort of combination of the approximation when $\lambda < c\mu$ and of the approximation when $\lambda > c\mu$ could make a nice choice for approximation of the distribution when $\lambda = c\mu$. In some sense, it is not surprising that these approximations are similar to the true limiting distribution, as the Erlang-A appears to be a $M/M/\infty$ queue with service rate $\mu$ (the approximation when $\lambda < c\mu$), when only considering the states up to $c$, and it also resembles some sort of shifted $M/M/\infty$ queue with service rate $\theta$ (which also describes the approximation when $\lambda \geq c\mu$) for states $c + 1$ and beyond.  Finally, it would be interesting to extend this to networks of Erlang-A queues like in \citet{pender2017approximating}, however, we would have to keep track of the routing probabilities carefully to keep track of the convexity/concavity of the rate functions.


\bibliographystyle{plainnat}
\bibliography{ErlangAJensens2}

\providecommand{\noopsort}[1]{} \providecommand{\doi}[1]{\texttt{doi:#1}}
  \providecommand{\available}[1]{Available at \texttt{#1}}
  \providecommand{\availablet}[2]{Available at \texttt{#2}}
\begin{thebibliography}{26}
\providecommand{\natexlab}[1]{#1}
\providecommand{\url}[1]{\texttt{#1}}
\expandafter\ifx\csname urlstyle\endcsname\relax
  \providecommand{\doi}[1]{doi: #1}\else
  \providecommand{\doi}{doi: \begingroup \urlstyle{rm}\Url}\fi

\bibitem[Eick et~al.(1993)Eick, Massey, and Whitt]{eick1993mt}
Stephen~G Eick, William~A Massey, and Ward Whitt.
\newblock {$M_t/G/\infty$} queues with sinusoidal arrival rates.
\newblock \emph{Management Science}, 39\penalty0 (2):\penalty0 241--252, 1993.

\bibitem[Engblom and Pender(2014)]{engblom2014approximations}
Stefan Engblom and Jamol Pender.
\newblock Approximations for the moments of nonstationary and state dependent
  birth-death queues.
\newblock \emph{arXiv preprint arXiv:1406.6164}, 2014.

\bibitem[Ferragut and Paganini(2012)]{ferragut2012content}
Andres Ferragut and Fernando Paganini.
\newblock Content dynamics in {P2P} networks from queueing and fluid
  perspectives.
\newblock In \emph{Proceedings of the 24th International Teletraffic Congress},
  page~11. International Teletraffic Congress, 2012.

\bibitem[Hale and Lunel(2013)]{hale2013introduction}
Jack~K Hale and Sjoerd M~Verduyn Lunel.
\newblock \emph{Introduction to functional differential equations}, volume~99.
\newblock Springer Science \& Business Media, 2013.

\bibitem[Halfin and Whitt(1981)]{halfin1981heavy}
Shlomo Halfin and Ward Whitt.
\newblock Heavy-traffic limits for queues with many exponential servers.
\newblock \emph{Operations Research}, 29\penalty0 (3):\penalty0 567--588, 1981.

\bibitem[Hampshire and Massey(2010)]{hampshire2010dynamic}
Robert~C Hampshire and William~A Massey.
\newblock Dynamic optimization with applications to dynamic rate queues.
\newblock In \emph{Risk and Optimization in an Uncertain World}, pages
  208--247. INFORMS, 2010.

\bibitem[Hampshire et~al.(2009{\natexlab{a}})Hampshire, Jennings, and
  Massey]{hampshire2009time}
Robert~C Hampshire, Otis~B Jennings, and William~A Massey.
\newblock A time-varying call center design via {L}agrangian mechanics.
\newblock \emph{Probability in the Engineering and Informational Sciences},
  23\penalty0 (02):\penalty0 231--259, 2009{\natexlab{a}}.

\bibitem[Hampshire et~al.(2009{\natexlab{b}})Hampshire, Massey, and
  Wang]{hampshire2009dynamic}
Robert~C Hampshire, William~A Massey, and Qiong Wang.
\newblock Dynamic pricing to control loss systems with quality of service
  targets.
\newblock \emph{Probability in the Engineering and Informational Sciences},
  23\penalty0 (02):\penalty0 357--383, 2009{\natexlab{b}}.

\bibitem[Ko and Pender(2016)]{ko2016strong}
Young~Myoung Ko and Jamol Pender.
\newblock Strong approximations for time-varying infinite-server queues with
  non-renewal arrival and service processes.
\newblock \emph{Stochastic Models}, 2016.

\bibitem[Ko and Pender(2017)]{ko2017diffusion}
Young~Myoung Ko and Jamol Pender.
\newblock Diffusion limits for the {$(MAP_t/Ph_t/\infty)$ $N$} queueing
  network.
\newblock \emph{Operations Research Letters}, 45\penalty0 (3):\penalty0
  248--253, 2017.

\bibitem[Mandelbaum et~al.(1998)Mandelbaum, Massey, and
  Reiman]{mandelbaum1998strong}
Avi Mandelbaum, William~A Massey, and Martin~I Reiman.
\newblock Strong approximations for {M}arkovian service networks.
\newblock \emph{Queueing Systems}, 30\penalty0 (1):\penalty0 149--201, 1998.

\bibitem[Mandelbaum et~al.(2002)Mandelbaum, Massey, Reiman, Stolyar, and
  Rider]{mandelbaum2002queue}
Avi Mandelbaum, William~A Massey, Martin~I Reiman, Alexander Stolyar, and Brian
  Rider.
\newblock Queue lengths and waiting times for multiserver queues with
  abandonment and retrials.
\newblock \emph{Telecommunication Systems}, 21\penalty0 (2):\penalty0 149--171,
  2002.

\bibitem[Massey(2002)]{massey2002analysis}
William~A Massey.
\newblock The analysis of queues with time-varying rates for telecommunication
  models.
\newblock \emph{Telecommunication Systems}, 21\penalty0 (2):\penalty0 173--204,
  2002.

\bibitem[Massey and Pender(2011)]{massey2011poster}
William~A Massey and Jamol Pender.
\newblock Poster: skewness variance approximation for dynamic rate multiserver
  queues with abandonment.
\newblock \emph{ACM SIGMETRICS Performance Evaluation Review}, 39\penalty0
  (2):\penalty0 74--74, 2011.

\bibitem[Massey and Pender(2013)]{massey2013gaussian}
William~A Massey and Jamol Pender.
\newblock Gaussian skewness approximation for dynamic rate multi-server queues
  with abandonment.
\newblock \emph{Queueing Systems}, 75\penalty0 (2-4):\penalty0 243--277, 2013.

\bibitem[Massey and Pender(2017)]{massey2017performance}
William~A Massey and Jamol Pender.
\newblock Performance and provisioning analysis for the dynamic rate {Erlang-A}
  queue.
\newblock 2017.

\bibitem[Matis and Feldman(2001)]{matis2001transient}
Timothy~I Matis and Richard~M Feldman.
\newblock Transient analysis of state-dependent queueing networks via cumulant
  functions.
\newblock \emph{Journal of Applied Probability}, 38\penalty0 (4):\penalty0
  841--859, 2001.

\bibitem[Niyirora and Pender(2016)]{niyirora2016optimal}
Jerome Niyirora and Jamol Pender.
\newblock Optimal staffing in nonstationary service centers with constraints.
\newblock \emph{Naval Research Logistics (NRL)}, 63\penalty0 (8):\penalty0
  615--630, 2016.

\bibitem[Pender(2014{\natexlab{a}})]{pender2014gram}
Jamol Pender.
\newblock Gram-{C}harlier expansion for time varying multiserver queues with
  abandonment.
\newblock \emph{SIAM Journal on Applied Mathematics}, 74\penalty0 (4):\penalty0
  1238--1265, 2014{\natexlab{a}}.

\bibitem[Pender(2014{\natexlab{b}})]{pender2014laguerre}
Jamol Pender.
\newblock Laguerre polynomial expansions for time varying multiserver queues
  with abandonment, 2014{\natexlab{b}}.

\bibitem[Pender(2016)]{pender2016sampling}
Jamol Pender.
\newblock Sampling the functional {K}olmogorov forward equations for
  nonstationary queueing networks.
\newblock \emph{INFORMS Journal on Computing}, 29\penalty0 (1):\penalty0 1--17,
  2016.

\bibitem[Pender and Ko(2017)]{pender2017approximations}
Jamol Pender and Young~Myoung Ko.
\newblock Approximations for the queue length distributions of time-varying
  many-server queues.
\newblock \emph{INFORMS Journal on Computing}, 29\penalty0 (4):\penalty0
  688--704, 2017.

\bibitem[Pender and Massey(2017)]{pender2017approximating}
Jamol Pender and William~A Massey.
\newblock Approximating and stabilizing dynamic rate {J}ackson networks with
  abandonment.
\newblock \emph{Probability in the Engineering and Informational Sciences},
  31\penalty0 (1):\penalty0 1--42, 2017.

\bibitem[Pender and Phung-Duc(2016)]{pender2016law}
Jamol Pender and Tuan Phung-Duc.
\newblock A law of large numbers for {$M/M/c/Delay$} off-setup queues with
  nonstationary arrivals.
\newblock In \emph{International Conference on Analytical and Stochastic
  Modeling Techniques and Applications}, pages 253--268. Springer, 2016.

\bibitem[Qin and Pender(2017)]{qin2017dynamic}
Ziyuan Qin and Jamol Pender.
\newblock Dynamic control for nonstationary queueing networks.
\newblock 2017.

\bibitem[Yom-Tov and Mandelbaum(2014)]{yom2014erlang}
Galit~B Yom-Tov and Avishai Mandelbaum.
\newblock Erlang-{R}: A time-varying queue with reentrant customers, in support
  of healthcare staffing.
\newblock \emph{Manufacturing \& Service Operations Management}, 16\penalty0
  (2):\penalty0 283--299, 2014.

\end{thebibliography}

\end{document}